%% file: main.tex
\documentclass[10pt]{article}

\input{packages/import}
\input{acronyms/import}
\input{symbols/import}

\input{style/import}

\input{format/import}

\title{
    Distributionally Robust Linear Quadratic\\
    Gaussian Regulator with Stationary Distributions
}
\author{
    Alain Sch\"obi,
    Nicolas Lanzetti,
    Florian D\"orfler,
    Raffaello D'Andrea,
    and
    Antonio Terpin%
}
\abstractText{\input{sections/abstract.txt}}

\newcommand{\email}[1]{\mbox{#1}}%
\affiliationText{
    This work was supported as a part of NCCR Automation, a National Centre of Competence in Research, funded by the Swiss National Science Foundation (grant number 51NF40\_225155). This research was conducted when all authors were affiliated with ETH Z\"urich.

    Alain Sch\"obi is with the
    Risk Analytics and Optimization Chair at EPFL, Switzerland (\email{alain.schoebi@epfl.ch}).
    Nicolas Lanzetti is with the Department of Computing and Mathematical Sciences at Caltech, USA (\email{lnicolas@caltech.edu}).
    Florian D\"orfler is with the Automatic Control Laboratory
    at ETH Z\"urich, Switzerland (\email{dorfler@ethz.edu}).
    Raffaello D'Andrea and Antonio Terpin are with
    the Institute for Dynamic Systems and Control at ETH Z\"urich, Switzerland (\email{\{rdandrea,aterpin\}@ethz.ch}).
}

\begin{document}

\maketitle
\makeabstract
\makeaffiliations

\input{sections/import}

\bibliography{references}

\end{document}

%% file: packages/import.tex
\usepackage{xcolor}

\usepackage{amsmath,amssymb,mathtools,amsthm,aligned-overset,physics}

\usepackage[hidelinks]{hyperref}

\usepackage[capitalize]{cleveref}
\usepackage{pifont}

\allowdisplaybreaks

\newtheorem{theorem}{Theorem}

\newtheorem{assumption}{Assumption}

\usepackage[framemethod=tikz]{mdframed}

\usepackage{pgfplots}
\pgfplotsset{compat=1.18}

\usepackage{dsfont}
\usepackage{mathrsfs}

\usepackage{xspace}
\usepackage{xstring}

\usepackage[normalem]{ulem}

\usepackage{algorithm}
\usepackage{algpseudocode}


\usepackage{enumitem}

%% file: acronyms/import.tex
\usepackage[acronym]{glossaries}
\newacronym{acr:lqg}{LQG}{Linear Quadratic Gaussian}
\newacronym{acr:lqr}{LQR}{Linear Quadratic Regulator}
\newacronym{acr:kf}{KF}{Kalman filter}
\newacronym{acr:sls}{SLS}{System Level Synthesis}
\newacronym{acr:fw}{FW}{Frank-Wolfe}
\newacronym{acr:dr}{DR}{distributionally robust}
\newacronym{acr:ibr}{IBR}{iterated best-response}

%% file: symbols/import.tex
\input{symbols/math}
\input{symbols/utils}

\input{symbols/text}
\input{symbols/system}
\input{symbols/policy}
\input{symbols/batch}
\input{symbols/dro}
\input{symbols/cost}
\input{symbols/algorithms}
\input{symbols/icons}

%% file: symbols/math.tex
\let\originalexists\exists
\renewcommand{\exists}{\originalexists\,}
\newcommand{\middleSt}{\,\middle|\,}

\renewcommand{\dim}{d}

\newcommand{\equalByDefRight}{=\mathrel{\mathop:}}

\newcommand{\openClosedInterval}[2]{\mathopen{(}#1,#2\mathclose{]}}

\newcommand{\normalSymbol}{\mathcal{N}}
\NewDocumentCommand{\normal}{s g g g}{%
  \ensuremath{%
    \IfValueTF{#2}{%
      \IfValueTF{#3}{%
        \IfValueTF{#4}{%
          \IfBooleanTF{#1}{%
            \normalSymbol\!\left(#2; #3, #4\right)%
          }{%
            \normalSymbol(#2; #3, #4)%
          }%
        }{%
          \IfBooleanTF{#1}{%
            \normalSymbol\!\left(#2, #3\right)%
          }{%
            \normalSymbol(#2, #3)%
          }%
        }%
      }{%
        \normalSymbol(?)%
      }%
    }{%
      \normalSymbol%
    }%
  }%
}

\newcommand{\expected}[2]{%
   \mathbb{E}_{\substack{#1}}\!%
   \left[%
      #2
   \right]%
}

\NewDocumentCommand{\probabilityMeasure}{g}{%
  \IfValueTF{#1}%
  {\mathbb{P}_{#1}}%
  {\mathbb{P}}%
}

\NewDocumentCommand{\probabilitySet}{sm}{
\ensuremath{
  \IfBooleanTF{#1}
    {\mathcal{P}(#2)}
    {\mathcal{P}\!\left(#2\right)}
}}

\NewDocumentCommand{\probabilitySetFiniteSecondMoment}{sm}{
\ensuremath{
  \IfBooleanTF{#1}
    {\mathcal{P}_{2}(#2)}
    {\mathcal{P}_{2}\!\left(#2\right)}
}}

\newcommand{\productMeasure}{\otimes}

\newcommand{\coupling}{\gamma}
\newcommand{\setPlans}[2]{\Gamma(#1,#2)}

\NewDocumentCommand{\wassersteinDistance}{s g g}{%
  \ensuremath{%
    \IfBooleanTF{#1}{%
      \IfValueTF{#2}{%
        \IfValueTF{#3}
          {\mathds{W}\!\left(#2, #3\right)}%
          {\mathds{W}\!\left(#2, ?\right)}%
      }{\mathds{W}}%
    }{%
      \IfValueTF{#2}{%
        \IfValueTF{#3}
          {\mathds{W}(#2, #3)}%
          {\mathds{W}(#2, ?)}%
      }{\mathds{W}}%
    }%
  }%
}

\NewDocumentCommand{\gelbrichDistance}{s g g}{%
  \ensuremath{%
    \IfBooleanTF{#1}{%
      \IfValueTF{#2}{%
        \IfValueTF{#3}
          {\mathds{G}\!\left(#2, #3\right)}%
          {\mathds{G}\!\left(#2, ?\right)}%
      }{\mathds{G}}%
    }{%
      \IfValueTF{#2}{%
        \IfValueTF{#3}
          {\mathds{G}(#2, #3)}%
          {\mathds{G}(#2, ?)}%
      }{\mathds{G}}%
    }%
  }%
}

\newcommand{\probabilityMeasureA}{\mathbb{P}_{1}}
\newcommand{\probabilityMeasureB}{\mathbb{P}_{2}}
\newcommand{\meanVector}{\mu}

\newcommand{\covarianceMatrix}{\Sigma}
\newcommand{\covarianceMatrixA}{\covarianceMatrix_{1}}
\newcommand{\covarianceMatrixB}{\covarianceMatrix_{2}}

\NewDocumentCommand{\generalMathOperatorCommand}{mgg}{
 \ensuremath{
    \IfValueTF{#2}{%
        \IfValueTF{#3}{%
            \operatorname*{#1}_{#2}\left\{{#3}\right\}
        }{%
            \operatorname*{#1}_{#2}
        }%
    }{%
        \operatorname*{#1}%
    }%
}}
\RenewDocumentCommand{\max}{gg}{\generalMathOperatorCommand{max}{#1}{#2}}
\RenewDocumentCommand{\min}{gg}{\generalMathOperatorCommand{min}{#1}{#2}}
\NewDocumentCommand{\argmax}{gg}{\generalMathOperatorCommand{argmax}{#1}{#2}}
\NewDocumentCommand{\argmin}{gg}{\generalMathOperatorCommand{argmin}{#1}{#2}}
\RenewDocumentCommand{\inf}{gg}{\generalMathOperatorCommand{inf}{#1}{#2}}
\RenewDocumentCommand{\sup}{gg}{\generalMathOperatorCommand{sup}{#1}{#2}}

\newcommand{\diag}[1]{\mathrm{diag}(#1)}
\newcommand{\mineigenvalue}[1]{\lambda_\mathrm{min}(#1)}
\newcommand{\maxeigenvalue}[1]{\lambda_\mathrm{max}(#1)}
\NewDocumentCommand{\identity}{g}{
    \IfValueTF{#1}
    {I_{#1}}
    {I}
}
\RenewDocumentCommand{\trace}{sm}{
\ensuremath{
  \IfBooleanTF{#1}
    {\mathrm{tr}\!\left(#2\right)}
    {\mathrm{tr}(#2)}
}}

\newcommand{\half}[1]{{#1}^{\frac{1}{2}}}

\NewDocumentCommand{\inverse}{sm}{
\ensuremath{
  \IfBooleanTF{#1}
    {\left(#2\right)^{-1}}
    {#2^{-1}}
}}

\NewDocumentCommand{\R}{g}{%
  \IfValueTF{#1}{%
    \IfSubStr{#1}{,}{
      \StrBefore{#1}{,}[\FirstDim]%
      \StrBehind{#1}{,}[\SecondDim]%
      \mathbb{R}^{\FirstDim \times \SecondDim}%
    }%
    {
      \mathbb{R}^{#1}
    }%
  }%
  {%
       \mathbb{R}%
  }%
}

\newcommand{\posSemiDef}[1]{\mathbb S_{+}^{#1}}

%% file: symbols/utils.tex
\NewDocumentCommand{\todo}{g}{%
  \IfValueTF{#1}{%
    \textcolor{Fuchsia}{\textbf{TODO}: #1}\xspace%
  }%
  {%
    \textcolor{Fuchsia}{\textbf{TODO}}\\\xspace%
  }%
}

%% file: symbols/text.tex
\newcommand{\ie}{i.e.,\xspace}

%% file: symbols/system.tex
\newcommand{\stateDim}{n}
\newcommand{\controlInputDim}{m}
\newcommand{\observationDim}{p}
\newcommand{\horizon}{T}
\newcommand{\processDim}{\stateDim}
\newcommand{\measurementDim}{\observationDim}

\renewcommand{\t}{t}

\newcommand{\tTau}{\tau}

\makeatletter

\newcommand{\inFromTo}[2]{=#1, \ldots, #2}
\newcommand{\inZeroTo}[1]{\inFromTo{0}{#1}}

\newcommand{\inT}{\inZeroTo{\horizon}}
\newcommand{\tInT}{\@ifstar{\t & \inT}{\t \inT}}
\newcommand{\tInTMinusOne}{\@ifstar{\t & \inZeroTo{\horizon-1}}{\t \inZeroTo{\horizon-1}}}

\makeatother

\newcommand{\state}[1]{x_{#1}}
\newcommand{\controlInput}[1]{u_{#1}}
\NewDocumentCommand{\observation}{g}{
    \IfValueTF{#1}{y_{#1}}%
    {y}%
}

\newcommand{\initialState}{\mathrm{x}}

\newcommand{\initialStateCovariance}{{\covarianceMatrix}_{\initialState}}
\newcommand{\initialStateCovarianceRef}{\widehat{\covarianceMatrix}_{\initialState}}

\newcommand{\initialStateDistribution}{\probabilityMeasure{\initialState}}
\newcommand{\initialStateDistributionRef}{\widehat{\probabilityMeasure}_{\initialState}}

\NewDocumentCommand{\processNoise}{g}{
    \IfValueTF{#1}{v_{#1}}%
    {v}%
}

\newcommand{\processNoiseTrajectory}{\boldsymbol{\processNoise{}}}

\NewDocumentCommand{\processMean}{g}{\ensuremath{{\meanVector}_{\processNoise{#1}}}}
\NewDocumentCommand{\processCovariance}{g}{\ensuremath{{\covarianceMatrix}_{\processNoise{#1}}}}
\NewDocumentCommand{\processCovarianceRef}{g}{\widehat{\covarianceMatrix}_{\processNoise{#1}}}
\NewDocumentCommand{\processCovarianceOptimal}{g}{\ensuremath{{\covarianceMatrix}_{\processNoise{#1}}^{\ast}}}

\NewDocumentCommand{\processDistribution}{g}{\probabilityMeasure{\processNoise{#1}}}
\NewDocumentCommand{\processDistributionRef}{g}{\widehat{\probabilityMeasure}_{\processNoise{#1}}}

\NewDocumentCommand{\measurementNoise}{g}{%
  \IfValueTF{#1}{w_{#1}}%
  {w}%
}
\NewDocumentCommand{\measurementNoiseRef}{g}{\widehat{\measurementNoise}_{#1}}
\newcommand{\measurementNoiseTrajectory}{\boldsymbol{\measurementNoise{}}}

\NewDocumentCommand{\measurementMean}{g}{\ensuremath{{\meanVector}_{\measurementNoise{#1}}}}
\NewDocumentCommand{\measurementCovariance}{g}{\ensuremath{{\covarianceMatrix}_{\measurementNoise{#1}}}}
\NewDocumentCommand{\measurementCovarianceRef}{g}{\widehat{\covarianceMatrix}_{\measurementNoise{#1}}}

\NewDocumentCommand{\measurementDistribution}{g}{\probabilityMeasure{\measurementNoise{#1}}}
\NewDocumentCommand{\measurementDistributionRef}{g}{\widehat{\probabilityMeasure}_{\measurementNoise{#1}}}

\newcommand{\allVariables}[2]{%
    #1{0}, \ldots, #1{#2}
}

\newcommand{\allProcessNoises}{\allVariables{\processNoise}{\horizon-1}}
\newcommand{\allMeasurementNoises}{\allVariables{\measurementNoise}{\horizon-1}}

\newcommand{\allDisturbancesCovariance}{
    \mathbf{\covarianceMatrix}
}
\newcommand{\allDisturbancesCovarianceRef}{
    \widehat{\mathbf{\covarianceMatrix}}
}
\newcommand{\allDisturbancesCovarianceStar}{
    \allDisturbancesCovariance^{\ast}
}

\newcommand{\jointDistribution}{\mathbb{P}}
\newcommand{\jointDistributionStar}{\jointDistribution^\ast}
\newcommand{\jointDistributionRef}{\widehat{\jointDistribution}}
\newcommand{\jointDistributionTrue}{\jointDistribution^\textrm{true}}

\newcommand{\jointDistributionFactorizationStationary}{%
\initialStateDistribution
\productMeasure
\left(\processDistribution{}\right)^{\productMeasure \horizon}
\productMeasure
\left(\measurementDistribution{}\right)^{\productMeasure \horizon}
}

\newcommand{\probabilityMeasureTemplate}{\probabilityMeasure{\templateVariable}}
\newcommand{\probabilityMeasureRefTemplate}{\widehat{\probabilityMeasure}_{\templateVariable}}

\newcommand{\Amat}[1]{A_{#1}}
\newcommand{\Bmat}[1]{B_{#1}}
\newcommand{\Cmat}[1]{C_{#1}}

\newcommand{\stageStateCostMatrix}[1]{Q_{#1}}
\newcommand{\terminalStateCostMatrix}{\stageStateCostMatrix{\horizon}}
\newcommand{\stageInputCostMatrix}[1]{R_{#1}}

\newcommand{\quadraticCostExpression}{%
    \state{\horizon}^\top \terminalStateCostMatrix\state{\horizon} + \sum_{\t=0}^{\horizon-1} \state{\t}^\top \stageStateCostMatrix{\t}\state{\t} + \controlInput{\t}^\top \stageInputCostMatrix{\t}\controlInput{\t}
}

\makeatletter
\newcommand{\linearDynamicsExpression}{%
  \@ifstar
    {\state{\t+1} = \Amat{\t}\state{\t} + \Bmat{\t}\controlInput{\t} + \processNoise{\t}}%
    {\state{\t+1} &= \Amat{\t}\state{\t} + \Bmat{\t}\controlInput{\t} + \processNoise{\t}}%
}
\makeatother

\makeatletter
\newcommand{\observationModelExpression}{%
  \@ifstar
    {\observation{\t} = \Cmat{\t}\state{\t} + \measurementNoise{\t}}%
    {\observation{\t} &= \Cmat{\t}\state{\t} + \measurementNoise{\t}}%
}
\makeatother

%% file: symbols/policy.tex
\newcommand{\affineSymbol}{\mathrm{aff}}
\newcommand{\linearSymbol}{\mathrm{lin}}

\NewDocumentCommand{\policy}{g}{%
\ensuremath{%
  \IfValueTF{#1}%
  {\pi_{#1}}%
  {\boldsymbol{\pi}}%
}}
\NewDocumentCommand{\policyStar}{g}{%
\ensuremath{%
  \IfValueTF{#1}%
  {\policy{#1}^\ast}%
  {\boldsymbol{\policy}^\ast}%
}}

\newcommand{\policyUnpacked}{
    \policy = (\policy{0}, \ldots, \policy{\horizon -1})
}

\newcommand{\policySet}{\Pi}
\newcommand{\policySetAffine}{\policySet_{\affineSymbol}}
\newcommand{\policySetLinear}{\policySet_{\linearSymbol}}

\newcommand{\controlInputGain}[2]{K_{#1,#2}}
\newcommand{\controlInputGainMatrix}{K}
\newcommand{\controlInputGainMatrixStar}{K^\ast}

\newcommand{\controlInputGainMatrixR}{\R{\controlInputDim\horizon, \observationDim\horizon}}
\newcommand{\controlInputGainMatrixInR}{\controlInputGainMatrix \in \controlInputGainMatrixR}

\newcommand{\controlInputGainMatrixSet}{\mathcal{K}}
\newcommand{\controlInputGainMatrixSetCompact}{\controlInputGainMatrixSet'}
\newcommand{\controlInputGainMatrixInU}{\controlInputGainMatrix \in \controlInputGainMatrixSet}
\newcommand{\controlInputGainMatrixInUCompact}{\controlInputGainMatrix \in \controlInputGainMatrixSetCompact}

\newcommand{\controlInputOffset}[1]{q_{#1}}

%% file: symbols/batch.tex
\newcommand{\matBatchH}{H}
\newcommand{\matBatchG}{G}
\newcommand{\matBatchL}{L}
\newcommand{\matBatchC}{C}
\newcommand{\matBatchQ}{Q}
\newcommand{\matBatchB}{B}
\newcommand{\matBatchR}{R}

\newcommand{\matBatchCostInitialStateCovariance}{
    P_{\initialState}
}
\newcommand{\matBatchCostProcessCovariance}{
    P_{\processNoise{}}
}
\newcommand{\matBatchCostMeasurementCovariance}{
    P_{\measurementNoise{}}
}
\newcommand{\matBatchCostTemplateCovariance}{P_{\templateVariable}}

\newcommand{\matBatchCostProcessCovarianceTrajectory}{
    \boldsymbol{P_{\processNoiseTrajectory}}
}
\newcommand{\matBatchCostMeasurementCovarianceTrajectory}{
    \boldsymbol{P_{\measurementNoiseTrajectory}}
}

\newcommand{\HU}{\matBatchH \controlInputGainMatrix}
\newcommand{\HUCPlusICompact}{(\matBatchH \controlInputGainMatrix \matBatchC \! + \! \identity)}
\newcommand{\UCL}{\controlInputGainMatrix \matBatchC \matBatchL}
\newcommand{\UCG}{\controlInputGainMatrix \matBatchC \matBatchG}

\newcommand{\matBatchCostInitialStateCovarianceExpressionCompact}{
    (\UCL)^{\!\top} \matBatchR (\UCL)
    +
    \matBatchL^{\!\top} \HUCPlusICompact^{\!\top}
    \matBatchQ
    \HUCPlusICompact \matBatchL
}

\newcommand{\matBatchCostProcessCovarianceExpressionCompact}{
    (\UCG)^{\!\top} \matBatchR (\UCG)
    +
    \matBatchG^{\!\top} \HUCPlusICompact^{\!\top}
    \matBatchQ
    \HUCPlusICompact \matBatchG
}

\newcommand{\matBatchCostMeasurementCovarianceExpressionCompact}
{
    \controlInputGainMatrix^{\!\top}
    \matBatchR
    \controlInputGainMatrix
    +
    (\HU)^{\!\top}
    \matBatchQ
    (\HU)
}

%% file: symbols/dro.tex
\NewDocumentCommand{\expectedCost}{gg}{%
\ensuremath{
  \IfValueTF{#1}{%
      \IfValueTF{#2}{%
          \mathcal{J}\!\left(#1, #2\right)%
      }%
      {%
          \mathcal{J}(?)%
      }%
  }%
  {%
      \mathcal{J}%
  }%
}}
\newcommand{\expectedCostPiP}{\expectedCost{\policy}{\jointDistribution}}

\newcommand{\dualTemplate}{\gamma_{\templateVariable}}

\newcommand{\radius}[1]{r_{#1}}

\NewDocumentCommand{\radiusProcess}{g}{\radius{\processNoise{#1}}}
\NewDocumentCommand{\radiusMeasurement}{g}{\radius{\measurementNoise{#1}}}

\newcommand{\unconstrainedMeanSymbol}{\odot}

\NewDocumentCommand{\ambiguitySet}{gg}{%
  \IfValueTF{#1}{%
    \IfValueTF{#2}
      {\mathcal{W}_{#1}^{#2}}%
      {\mathcal{W}_{#1}}}%
  {\mathcal{W}}%
}

\newcommand{\ambiguitySetTemplate}{\ambiguitySet{\templateVariable}}
\newcommand{\ambiguitySetInitialState}{\ambiguitySet{\initialState}}
\NewDocumentCommand{\ambiguitySetProcess}{g}{
  \IfValueTF{#1}{%
    \ambiguitySet{\processNoise{#1}}%
  }{%
    \ambiguitySet{\processNoise}%
  }%
}
\NewDocumentCommand{\ambiguitySetMeasurement}{g}{
  \IfValueTF{#1}{%
    \ambiguitySet{\measurementNoise{#1}}%
  }{%
    \ambiguitySet{\measurementNoise}%
  }%
}
\newcommand{\ambiguitySetGaussian}{\ambiguitySet{\normal}{}}

\newcommand{\ambiguitySetUnconstrainedMean}{\ambiguitySet{}{\unconstrainedMeanSymbol}}
\newcommand{\ambiguitySetGaussianUnconstrainedMean}{\ambiguitySet{\normal}{\unconstrainedMeanSymbol}}

\NewDocumentCommand{\ambiguitySetGelbrich}{gg}{%
  \IfValueTF{#1}{%
    \IfValueTF{#2}
      {\mathcal{G}_{#1}^{#2}}%
      {\mathcal{G}_{#1}}}%
  {\mathcal{G}}%
}

\newcommand{\ambiguitySetGelbrichCovariance}{\ambiguitySetGelbrich{}}
\newcommand{\ambiguitySetGelbrichCovariancePosDef}{\ambiguitySetGelbrich{}'}

\newcommand{\ambiguitySetGelbrichInitialState}{\ambiguitySetGelbrich{\initialState}}
\newcommand{\ambiguitySetGelbrichProcess}{\ambiguitySetGelbrich{\processNoise}}
\newcommand{\ambiguitySetGelbrichMeasurement}{\ambiguitySetGelbrich{\measurementNoise}}
\newcommand{\ambiguitySetGelbrichTemplate}{\ambiguitySetGelbrich{\templateVariable}}

\newcommand{\ambiguitySetGelbrichInitialStatePosDef}{\ambiguitySetGelbrichInitialState'}
\newcommand{\ambiguitySetGelbrichProcessPosDef}{\ambiguitySetGelbrichProcess'}
\newcommand{\ambiguitySetGelbrichMeasurementPosDef}{\ambiguitySetGelbrichMeasurement'}
\newcommand{\ambiguitySetGelbrichTemplatePosDef}{\ambiguitySetGelbrichTemplate'}

%% file: symbols/cost.tex
\NewDocumentCommand{\expectedCostLinearStationary}{gg}{%
\ensuremath{
    \IfValueTF{#1}{%
        \IfValueTF{#2}{%
            \expectedCost_{\linearSymbol}\!\left(#1, #2\right)
        }%
        {%
            \expectedCost_{\linearSymbol}(?)%
        }%
    }%
    {%
        \expectedCost_{\linearSymbol}%
    }%
}}

\newcommand{\expectedCostLinearStationaryUSigma}{
    \expectedCostLinearStationary{
        \controlInputGainMatrix
    }
    {
        \allDisturbancesCovariance
    }
}

\newcommand{\expectedCostLinearStationaryUSigmaExpression}{
    \trace{\matBatchCostInitialStateCovariance \initialStateCovariance}
    +
    \trace{\matBatchCostProcessCovariance \processCovariance}
    +
    \trace{\matBatchCostMeasurementCovariance \measurementCovariance}
}

%% file: symbols/algorithms.tex
\newcommand{\iterateNumber}[2]{#2^{(#1)}}

\newcommand{\stepSizeEngineer}[1]{\alpha_{#1}}
\newcommand{\stepSizeDevil}[1]{\beta_{#1}}

\newcommand{\delay}{\delta}

\newcommand{\controlInputGainMatrixIterate}[1]{\iterateNumber{#1}{\controlInputGainMatrix}}
\newcommand{\allDisturbancesCovarianceIterate}[1]{\iterateNumber{#1}{\allDisturbancesCovariance}}

\NewDocumentCommand{\minCostFunction}{g}{%
\ensuremath{
  \IfValueTF{#1}{\mathcal{F}(#1)}%
  {\mathcal{F}}%
}}
\newcommand{\minCostFunctionSigmas}{
    \minCostFunction{\allDisturbancesCovariance}
}

\newcommand{\BRMaxSymbol}{\textrm{max}}
\newcommand{\BRMinSymbol}{\textrm{min}}

\newcommand{\templateVariable}{z}
\newcommand{\templateVariableInXVW}{\templateVariable \in \{\initialState, \processNoise, \measurementNoise\}}
\newcommand{\covarianceTemplate}{{\covarianceMatrix}_\templateVariable}
\newcommand{\covarianceTemplateRef}{\widehat{\covarianceMatrix}_{\templateVariable}}
\newcommand{\radiusTemplate}{\radius{\templateVariable}}

\newcommand{\covarianceTemplateBR}{{\covarianceMatrix}_\templateVariable^{\BRMaxSymbol}}
\newcommand{\allDisturbancesCovarianceBR}{\allDisturbancesCovariance^{\BRMaxSymbol}}
\newcommand{\controlInputGainMatrixBR}{\controlInputGainMatrix^{\BRMinSymbol}}

%% file: symbols/icons.tex
\NewDocumentCommand{\engineerSymbol}{s}{%
  \IfBooleanTF{#1}%
    {%
        \raisebox{-0.2em}{\includegraphics[width=1.4em]{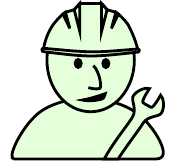}}\xspace%
    }%
    {%
        \smash{\raisebox{-0.2em}{\includegraphics[width=1.4em]{figures/engineer.pdf}}}\xspace%
    }%
}

\NewDocumentCommand{\devilSymbol}{s}{%
  \IfBooleanTF{#1}%
    {%
        \raisebox{-0.2em}{\includegraphics[width=1.4em]{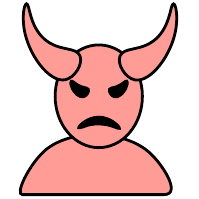}}\xspace%
    }%
    {%
        \smash{\raisebox{-0.2em}{\includegraphics[width=1.4em]{figures/devil.pdf}}}\xspace%
    }%
}

%% file: style/import.tex
\input{style/boxes}
\input{style/colors}

%% file: style/boxes.tex
\newenvironment{boxedassumption}
    {%
    \begin{mdframed}[
        innertopmargin=5pt,
        innerbottommargin=3pt,
        hidealllines=true,
        nobreak=true,
        backgroundcolor=gray!10,
    ]%
    \vspace{-\topsep}
    \begin{assumption}\strut%
    }
    {%
    \end{assumption}%
    \end{mdframed}%
}

\newenvironment{boxedtheorem}[1][]
    {%
    \begin{mdframed}[
        innertopmargin=5pt,
        innerbottommargin=3pt,
        hidealllines=true,
        nobreak=true,
        backgroundcolor=cost2!9
    ]%
    \vspace{-\topsep}
    \setlist[enumerate]{
      topsep=2pt,
      partopsep=0pt,
      itemsep=2pt,
      parsep=0pt,
    }%
    \begin{theorem}[#1]\hspace{-0.4em}\strut%
    }
    {%
    \end{theorem}%
    \end{mdframed}%
}

\newenvironment{boxedcontributions}
    {%
    \begin{mdframed}[
        innertopmargin=5pt,
        innerbottommargin=3pt,
        hidealllines=true,
        nobreak=true,
        backgroundcolor=robust!8
    ]%
    \noindent\textbf{Contributions. }\itshape\strut%
    }
    {%
    \end{mdframed}%
}

%% file: style/colors.tex
\definecolor{physical}{RGB}{112,112,112}
\definecolor{reference}{RGB}{245,151,9}
\definecolor{robust}{RGB}{135,20,202}

\definecolor{process}{RGB}{0,51,204}
\definecolor{measurement}{RGB}{0,204,51}

\definecolor{engineer}{HTML}{E5FFDA}
\definecolor{devil}{HTML}{FF9C95}

\colorlet{engineerHighlight}{engineer!60}
\colorlet{devilHighlight}{devil!30}

\definecolor{cost}{HTML}{0F4690}
\colorlet{costHighlight}{cost!20}

\definecolor{cost2}{HTML}{31A2D8}

\colorlet{region}{cost}

\newcommand{\colorboxEngineer}[1]{%
    {%
        \setlength{\fboxsep}{3pt}%
        \smash{\strut%
        \colorbox{engineerHighlight}{#1\vphantom{x}}}\xspace%
    }%
}
\newcommand{\colorboxDevil}[1]{%
    {%
        \setlength{\fboxsep}{3pt}%
        \smash{\strut%
        \colorbox{devilHighlight}{#1\vphantom{x}}}\xspace%
    }%
}
\newcommand{\colorboxSmall}[2]{%
    {%
        \setlength{\fboxsep}{1.5pt}%
        \smash{\strut%
        \colorbox{#1}{#2\vphantom{x}}}\xspace%
    }%
}

\newcommand{\colorboxEngineerEquation}[1]{
    {%
    \setlength{\fboxsep}{11pt}%
        \colorbox{engineerHighlight}{%
            $\displaystyle%
            \hspace{-6pt}%
                \smash{%
                    #1%
                }%
            \hspace{-6pt}%
            $%
        }%
    }
}

\newcommand{\colorboxDevilEquation}[1]{
    {%
    \setlength{\fboxsep}{11pt}%
        \colorbox{devilHighlight}{%
            $\displaystyle%
            \hspace{-6pt}%
                \smash{%
                    #1%
                }%
            \hspace{-6pt}%
            $%
        }%
    }
}

%% file: format/import.tex
\usepackage[
  letterpaper,
  textwidth=43pc,
  textheight=58pc,
  top=0.60in,
  heightrounded
]{geometry}

\makeatletter
\g@addto@macro\normalsize{%
  \setlength{\abovedisplayskip}{5pt plus 1pt minus 1pt}%
  \setlength{\belowdisplayskip}{5pt plus 1pt minus 1pt}%
  \setlength{\abovedisplayshortskip}{3pt plus 1pt minus 1pt}%
  \setlength{\belowdisplayshortskip}{4pt plus 1pt minus 1pt}%
}
\makeatother

\usepackage[final]{microtype}
\AtBeginDocument{%
  \sloppy
  \flushbottom
}

\usepackage[T1]{fontenc}

\usepackage{cite}
\usepackage{indentfirst}

\usepackage{lettrine}
\newcommand{\IEEEPARstart}[2]{%
  \lettrine[
    lines=2,
    loversize=0.10,
    lhang=0,
    nindent=0pt,
    findent=0.10em
  ]{\bfseries #1}{\MakeUppercase{#2}}%
}

\usepackage[title]{appendix}
\let\oldappendix\appendix

\renewcommand{\appendix}{
  \oldappendix
  \titleformat{\section}
    {\centering\sffamily\scshape\color{black}\fontsize{10}{12}\selectfont}
    {\appendixname~\thesection.}       
    {0.75em}
    {}
}

\usepackage[font=footnotesize]{caption}
\usepackage{titling}
\pretitle{%
  \begin{center}
  \fontsize{24}{28}\selectfont
}
\posttitle{%
  \par
  \end{center}
  \vspace{0.6em}
}

\preauthor{%
  \begin{center}
  \fontsize{12}{1}\selectfont
}
\postauthor{%
  \par
  \end{center}
  \vspace{0.4cm}
}

\predate{}
\postdate{}
\date{}

\let\oldmaketitle\maketitle
\renewcommand{\maketitle}{%
  \twocolumn[
    \oldmaketitle
  ]
}

\usepackage{titlesec}
\renewcommand{\thesection}{\Roman{section}}
\renewcommand{\thesubsection}{\Alph{subsection}}
\makeatletter
\def\p@subsection{\thesection-}
\makeatother

\titleformat{\section}
  {\centering\sffamily\scshape\color{black}\fontsize{10}{12}\selectfont}
  {\thesection.}
  {0.75em}
  {}

\titlespacing*{\section}
  {0pt}
  {3.5ex plus 0.8ex minus 0.5ex}
  {0.7ex plus 1ex}

\titleformat{\subsection}
  {\sffamily\itshape\color{black}\normalsize}
  {\thesubsection.}
  {0.75em}
  {}

\titlespacing*{\subsection}
  {0pt}
  {3.0ex plus 0.8ex minus 0.5ex}
  {0.7ex plus .5ex}

\newcommand{\makeabstract}{%
  Use \texttt{abstractText} to define the abstract.
}

\newcommand{\abstractText}[1]{%
    \renewcommand{\makeabstract}{
        \par\begingroup
        \sffamily\bfseries
        \fontsize{9}{10}\selectfont
        \setlength{\parindent}{0pt}%
        \setlength{\parskip}{0pt}%
        \noindent{\itshape Abstract}---\,#1\par
        \endgroup
    }
}

\newcommand{\makeaffiliations}{
    Use \texttt{affiliationText} first.
}

\newcommand{\affiliationText}[1]{
    \renewcommand{\makeaffiliations}{
        \begingroup
        \renewcommand{\thefootnote}{}%
        \footnotetext{%
            \hspace{-1.8em}%
            \begin{minipage}{\columnwidth}
            \footnotesize
            \setlength{\parindent}{1em}%
            \setlength{\parskip}{0pt}%
            #1
            \end{minipage}%
        }%
        \endgroup
    }
}

\renewenvironment{thebibliography}[1]{%
    \begin{oldthebibliography}{#1}%
    \fontsize{8pt}{9pt}\selectfont
    \setlength{\itemsep}{0pt}%
    \setlength{\parskip}{0pt}%
    \setlength{\parsep}{0pt}
}{%
    \end{oldthebibliography}%
}

%% file: sections/abstract.txt
We study the Linear Quadratic Gaussian regulation problem in the face of worst-case noise distributions when these are mutually independent, zero-mean, stationary, and within a radius (as measured by the Wasserstein distance) of some reference Gaussian noise distributions.
Compared with nonstationary ambiguity models, our stationary modeling choice reduces conservatism when the disturbances are stationary, as is often the case in practice.
We first show optimality of linear output-feedback policies via the existence of a Nash equilibrium in an equivalent zero-sum game between a control engineer and a fictitious adversary that compete to minimize and maximize the control cost.
Additionally, we show that Nash equilibria may fail to exist when the distributions are not restricted to have zero mean.
We then propose an iterated best-response algorithm to compute Nash equilibria and thereby the optimal feedback policies.
Our computational framework unifies two seemingly different viewpoints: game-theoretic iterated best response and Frank-Wolfe approaches.
Finally, we illustrate the robustness of the proposed controller to unknown noise distributions on an inverted pendulum with physical parameters, a canonical control benchmark.

%% file: sections/import.tex

\input{sections/introduction/import}

\input{sections/theory/import}

\input{sections/computation/import}

\input{sections/conclusion}

\input{sections/appendix/import}

%% file: sections/introduction/import.tex

\section{Introduction}
\label{sec:introduction}

\input{sections/introduction/setting}

\input{sections/introduction/disturbances}

\input{sections/introduction/numerical}

\input{sections/introduction/related}

%% file: sections/introduction/setting.tex

\IEEEPARstart{T}{he}
focus of this work is the optimal regulation of the (time-varying) discrete-time linear dynamical system
\begin{alignat}{2}
\label{eq:system-dynamics}
\linearDynamicsExpression, &\quad \tInTMinusOne*,
\intertext{
with state $\state{\t} \in \R{\stateDim}$, control input $\controlInput{\t} \in \R{\controlInputDim}$, process noise $\processNoise{\t} \in \R{\processDim}$, and system matrices $\Amat{\t}$, $\Bmat{\t}$.
We consider the case of imperfect state measurements via the observations $\observation{\t} \in \R{\observationDim}$,
}
\label{eq:observations}
\observationModelExpression, &\quad \tInTMinusOne*,
\end{alignat}
with measurement noise $\measurementNoise{\t} \in \R{\measurementDim}$, and observation matrix $\Cmat{\t}$.
Here, ``optimal'' regulation refers to minimizing the standard expected quadratic cost
$\expectedCostPiP$
defined as
\begin{equation}
\label{eq:the-cost}
    \expected
    {
        (\initialState, \processNoiseTrajectory, \measurementNoiseTrajectory) \sim \jointDistribution
    }
    {
        \quadraticCostExpression
    }
    ,
\end{equation}
where \eqref{eq:system-dynamics}
and
\eqref{eq:observations} hold,
and
$\state{0} = \initialState$,
$\controlInput{\t} = \policy{\t}(\observation{0}, \ldots, \observation{\t})$.

In \eqref{eq:the-cost}, $\stageStateCostMatrix{\t} \succeq 0$ and $\stageInputCostMatrix{\t} \succ 0$ are positive (semi-)definite cost matrices,
$\policyUnpacked \in \policySet$ is the feedback policy, which is a collection of feedback functions $\policy{\t}$ of all past observations $\observation{0}, \ldots, \observation{\t}$
(accordingly, $\policySet$ is the set of all feedback policies),
and $\jointDistribution$ is a \emph{joint} probability distribution over the \emph{exogenous random variables} describing the initial state $\initialState$, the process noises
$\processNoiseTrajectory = (\allProcessNoises)$,
and the measurement noises
$\measurementNoiseTrajectory = (\allMeasurementNoises)$.

\begin{figure}[t]
    \centering
    \vspace{-0.3cm}
    \input{figures/fig1}
    \caption{%
        Our distributionally robust \acrshort*{acr:lqg} control problem can be interpreted as a Stackelberg game between a control engineer (\engineerSymbol), trying to minimize the \colorboxSmall{costHighlight}{control cost}, and an adversarial Devil (\devilSymbol), seeking to maximize it.
        The game is sequential: the engineer first chooses their \colorboxSmall{engineerHighlight}{policy}; the Devil then selects the \colorboxSmall{devilHighlight}{noise distributions}.
        Under mild assumptions, we show that (i)~it is optimal for the engineer to select a linear policy, (ii)~the Devil reacts with Gaussian adversarial probability distributions, and (iii)~the game has a Nash~equilibrium.
    }
    \label{fig:motivating-example}
    \vspace{-.5cm}
\end{figure}

In practice, the exact noise distribution $\jointDistribution$ is unknown, and must be inferred from limited data or prior knowledge.
Rather than committing to a single nominal distribution, we consider
an \emph{ambiguity set}~$\ambiguitySet$ of plausible distributions~\cite{kuhn2025distributionally}.
We then seek a policy that \emph{%
    \colorboxEngineer{minimizes} the \colorboxDevil{worst-case} expected cost%
}
\begin{equation}
    \label{eq:dr-control-problem}
    \colorboxEngineerEquation{\inf{\policy \in \policySet}}
    \colorboxDevilEquation{\sup{\jointDistribution \in \ambiguitySet}}
    \,
    \expectedCostPiP
    .
\end{equation}

This \emph{\gls*{acr:dr}} control problem can be interpreted as a zero-sum Stackelberg game between a control engineer~(\engineerSymbol) and a fictitious adversary, the Devil (\devilSymbol), that compete to minimize/maximize the control cost \eqref{eq:dr-control-problem}.
The game runs sequentially: the control engineer first selects the \colorboxEngineer{policy $\policy$}; the Devil then chooses the
\colorboxDevil{noise distribution $\jointDistribution$} from the ambiguity set $\ambiguitySet$.
This is illustrated in~\cref{fig:motivating-example}.

%% file: figures/fig1.tex
\tikzstyle{myarrow}=[-,thick]
\tikzstyle{leaf}=[rectangle,draw=none,rounded corners]

\def\deltay{1.48}
\def\deltax{1.1}

\newcommand\drawBimodal{
\begin{tikzpicture}[rounded corners=0pt]
\begin{axis}[
        anchor=center,
        width=2.8cm,
        height=2cm,
        axis lines=none,
        domain=-3:3,
        xmin=-2.9,
        xmax=2.9,
        samples=100
    ]
    \addplot[thick] {0.5*exp(-(x+1.5)^2) + 0.5*exp(-(x-1.5)^2)};
\end{axis}
\end{tikzpicture}
}

\newcommand\drawGaussian{
\begin{tikzpicture}[rounded corners=0pt]
\begin{axis}[
        anchor=center,
        width=2.8cm,
        height=2cm,
        axis lines=none,
        ymin=-0.1,
        ymax=1.1,
        xmin=-2.9,
        xmax=2.9,
        domain=-3:3,
        samples=100
    ]
    \addplot[thick] {exp(-x^2)};
\end{axis}
\end{tikzpicture}
}
\pgfdeclarelayer{background}
\pgfsetlayers{background,main}

\begin{tikzpicture}
    \node[anchor=south] at (0,.5) {\includegraphics[width=.8\linewidth]{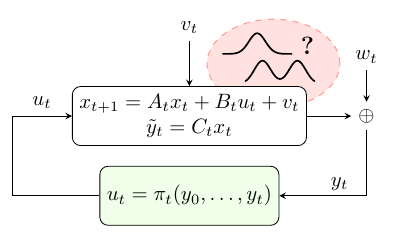}};
    \node[minimum size=0.1cm,inner sep=0pt,outer sep=0pt] (root) at (0,0) {\includegraphics[width=3em]{figures/engineer.pdf}};

    \node[minimum size=0.7cm,inner sep=0pt,outer sep=0pt] at (-2*\deltax,-\deltay) (n11) {\includegraphics[width=3em]{figures/devil.pdf}};
    \node[minimum size=0.7cm,inner sep=0pt,outer sep=0pt] at (+2*\deltax,-\deltay) (n12) {\includegraphics[width=3em]{figures/devil.pdf}};

    \node[leaf,fill=cost!45,rounded corners=0pt] at (-3*\deltax,-2*\deltay) (n21) {cost 3.5};
    \node[leaf,fill=cost!17,rounded corners=0pt] at (-\deltax,-2*\deltay) (n22) {cost 1.7};
    \node[leaf,fill=cost!58,rounded corners=0pt] at (\deltax,-2*\deltay) (n23) {cost 3.7};
    \node[leaf,fill=cost!9,rounded corners=0pt] at (3*\deltax,-2*\deltay) (n24) {cost 1.4};

    \node[circle] at (0,-0.8*\deltay) {\ldots};
    \node[circle] at (+2*\deltax,-1.8*\deltay) {\ldots};
    \node[circle] at (-2*\deltax,-1.8*\deltay) {\ldots};

    \coordinate (rshift) at ([yshift=-4pt] root);
    \coordinate (n11shift) at ([yshift=2pt] n11);
    \coordinate (n12shift) at ([yshift=2pt] n12);

    \begin{pgfonlayer}{background}
    \draw[myarrow] (rshift) -- (n11shift)
    node[pos=0.3,left,align=center,yshift=0.45cm,xshift=-0.45cm,fill=engineerHighlight,rounded corners]
    {
        Linear feedback\\[0.1cm]
        $\controlInput{t} = \mathbf{K}_t\!\left[\begin{smallmatrix}
            \observation{0} \\[-.2cm] \vdots \\ \observation{t}
        \end{smallmatrix}\right]$
    };

    \draw[myarrow] (rshift) -- (n12shift)
    node[pos=0.3,right,align=center,yshift=0.4cm,xshift=0.25cm,fill=engineerHighlight,rounded corners]
    {
        Nonlinear feedback \\ \small $\controlInput{t}=\policy{t}(\observation{0},\ldots,\observation{t})$
    };

    \draw[myarrow,shorten <=-8pt,shorten >=-2pt] (n11) -- (n21)
    node[pos=0.5,left,yshift=0.25cm,fill=devilHighlight,inner sep=2pt,outer sep=2pt,rounded corners]
    {
    \drawGaussian{}%
    };
    \draw[myarrow,shorten <=-8pt,shorten >=-2pt] (n11) -- (n22)
    node[pos=0.5,right,yshift=0.25cm,fill=devilHighlight,inner sep=2pt,outer sep=2pt,rounded corners]
    {
    \drawBimodal{}%
    };
    \draw[myarrow,shorten <=-8pt,shorten >=-2pt] (n12) -- (n23)
    node[pos=0.5,left,yshift=0.25cm,fill=devilHighlight,inner sep=2pt,outer sep=2pt,rounded corners]
    {
    \drawGaussian{}%
    };
    \draw[myarrow,shorten <=-8pt,shorten >=-2pt] (n12) -- (n24)
    node[pos=0.5,right,yshift=0.25cm,fill=devilHighlight,inner sep=2pt,outer sep=2pt,rounded corners]
    {
    \drawBimodal{}%
    };
    \end{pgfonlayer}

\end{tikzpicture}

%% file: sections/introduction/disturbances.tex

We make the following assumption on the disturbances.
\vspace{-0.05cm}
\begin{boxedassumption}
\label{ass:disturbances}
The initial state $\initialState$, the process noises $\allProcessNoises$, the measurement noises $\allMeasurementNoises$ are \emph{mutually independent} \emph{zero-mean} random variables.
Furthermore, the process and measurement noise sequences are \emph{stationary},
\ie $\processNoise{\t} \overset{\text{\tiny{i.i.d.}}}{\sim} \processDistribution$
and
$\measurementNoise{\t} \overset{\text{\tiny{i.i.d.}}}{\sim} \measurementDistribution$.
\end{boxedassumption}
\vspace{-0.05cm}
Mutual independence and stationarity
enable the factorization of the joint distribution $\jointDistribution$ into the product of marginals
\begin{equation*}
\jointDistribution
=
\jointDistributionFactorizationStationary
,
\end{equation*}
where $\cdot^{\productMeasure \horizon}$ denotes the $\horizon$-fold product measure.

To construct the ambiguity set, we consider a reference distribution $\jointDistributionRef$, which is typically inferred from available data.
We make the following assumption on~$\jointDistributionRef$.
\vspace{-0.05cm}
\begin{boxedassumption}
\label{ass:references}
    The reference $\smash{\jointDistributionRef}$ is taken as a \mbox{zero-mean} \emph{Gaussian} distribution, with marginals
    $\initialStateDistributionRef$, $\processDistributionRef$, $\measurementDistributionRef$,
    and
    corresponding covariance matrices
    $\initialStateCovarianceRef, \processCovarianceRef, \measurementCovarianceRef \succ 0$.
\end{boxedassumption}
\vspace{-0.05cm}
We then define the ambiguity set $\ambiguitySet$ as the collection of noise distributions $\jointDistribution$ that lie within a radius from $\jointDistributionRef$.
The distance between distributions is measured via the \emph{Wasserstein distance},
defined for two probability measures
$
\probabilityMeasureA, \probabilityMeasureB \in \probabilitySetFiniteSecondMoment{\R{\dim}}
$
as
\begin{equation*}
    \wassersteinDistance{\probabilityMeasureA}{\probabilityMeasureB}
    =
    \operatorname*{min}_{\coupling \in \setPlans{\probabilityMeasureA}{\probabilityMeasureB}}
    \left(
        \int_{\R{\dim}\times\R{\dim}} \norm{x - y}^2 \mathrm{d}\coupling(x, y)
    \right)
    ^{1/2},
\end{equation*}
where $\setPlans{\probabilityMeasureA}{\probabilityMeasureB}$ denotes the set of couplings between~$\probabilityMeasureA$ and~$\probabilityMeasureB$, \ie probability measures with marginals~$\probabilityMeasureA$ and~$\probabilityMeasureB$.

Formally, the ambiguity set $\ambiguitySet$ factors as Wasserstein balls
\begin{equation*}
    \ambiguitySet
    =
    \ambiguitySetInitialState
    \productMeasure
    (\ambiguitySetProcess)^{\productMeasure \horizon}
    \productMeasure
    (\ambiguitySetMeasurement)^{\productMeasure \horizon},
\end{equation*}
where, for $\templateVariableInXVW$ and ambiguity radius $\radiusTemplate > 0$,
\begin{equation*}
    \ambiguitySetTemplate
    =
    \left\{
        \probabilityMeasureTemplate
        \in \probabilitySetFiniteSecondMoment*{{\R}^{\textrm{dim}(\templateVariable)}}
        \middleSt
        \expected{
            \templateVariable \sim \probabilityMeasureTemplate
        }{\templateVariable} = 0
        ,\,
        \wassersteinDistance
        {\probabilityMeasureTemplate}
        {\probabilityMeasureRefTemplate}
        \le
        \radiusTemplate
    \right\}
    .
\end{equation*}

\begin{figure}[ht]
    \centering
    \includegraphics[width=0.31\textwidth]{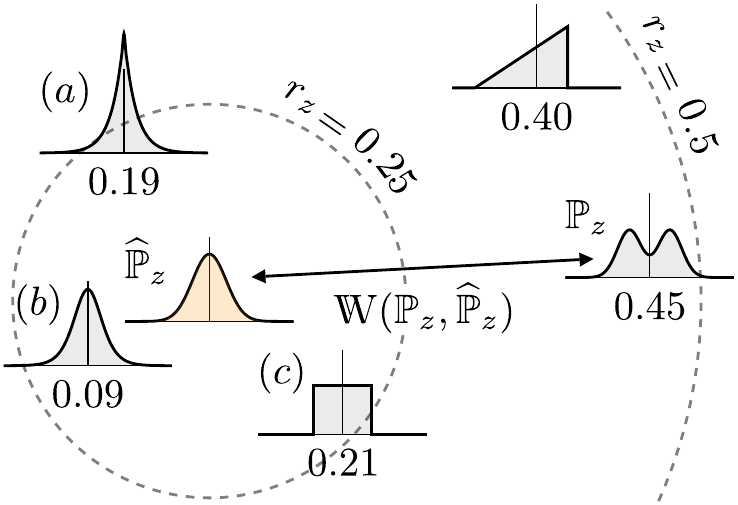}
    \vspace{-0.5\baselineskip}
    \caption{
        One-dimensional example of a Wasserstein ambiguity ball~$\ambiguitySet{\templateVariable}$ centered at the Gaussian \colorboxSmall{reference!20}{reference} distribution $\probabilityMeasureRefTemplate = \normal{0}{1}$.
        Candidate distributions $\probabilityMeasureTemplate$ displayed further from the reference correspond to larger Wasserstein distances, with the exact distance indicated below the PDF.
        The Laplace $(a)$, logistic $(b)$, and uniform $(c)$ distributions, all with zero mean and unit variance, lie within a ball of radius $\radiusTemplate = 0.25$.
    }
    \vspace{-.55cm}
    \label{fig:ambiguity-set}
\end{figure}

Importantly, as shown in \cref{fig:ambiguity-set}, $\ambiguitySet$ also contains non-Gaussian distributions, providing robustness beyond Gaussian disturbances.
In practice, each reference Gaussian $\probabilityMeasureRefTemplate$ can be calibrated from data and the radius $\radiusTemplate$ can be chosen using statistical concentration inequalities \cite{kuhn2025distributionally}.

\vspace{-0.05cm}
\begin{boxedcontributions}
We follow the minimax proof technique of~\cite{taskesen2024distributionally} to show that linear output-feedback policies remain optimal for the \gls*{acr:dr} control problem~\eqref{eq:dr-control-problem} with stationary disturbances.
The proof relies on the existence of a Nash equilibrium for the equivalent zero-sum game.
We further show that such an equilibrium may fail to exist when distributions are not zero-mean.
Then, we develop an intuitive algorithm for computing Nash equilibria and optimal policies.
The convergence analysis of our method unifies
two popular approaches in \gls*{acr:dr} control: iterated best response and Frank-Wolfe.
Finally, we demonstrate robustness to unknown noise distributions via simulations of an inverted pendulum with physical parameters.
\end{boxedcontributions}

\begin{figure}[ht]
    \centering
    \vspace{0.3em}
    \begin{tikzpicture}
    \draw[draw=none,fill=region,fill opacity=0.10]
      (-3.2,-0.95) rectangle (-2.6,1.0)
      node[region,xshift=-.27cm,yshift=-.2cm,text opacity=1] {\footnotesize R1};
    \draw[draw=none,fill=region,fill opacity=0.10]
      (-1.42,-0.95) rectangle (-.19,-0.1)
      node[region,xshift=-.6cm,yshift=-.2cm,text opacity=1] {\footnotesize R2};
    \node[anchor=base] at (-1.96,1.32) {\textbf{\textcolor{physical}{ideal}}};
    \node[anchor=base] at (0.065,1.32) {\textbf{\textcolor{reference}{nominal}}};
    \node[anchor=base] at (2.19,1.32) {\textbf{\textcolor{robust}{robust}}};
    \node[anchor=base] at (-1.95,-1.45) {expected cost};
    \node[anchor=base] at (2.28,-1.45) {average balancing time $[s]$};
    \node[rotate=90] at (-4.1,0.08) {frequency};
    \node[rectangle] at (0,0) {\includegraphics[width=0.47\textwidth]{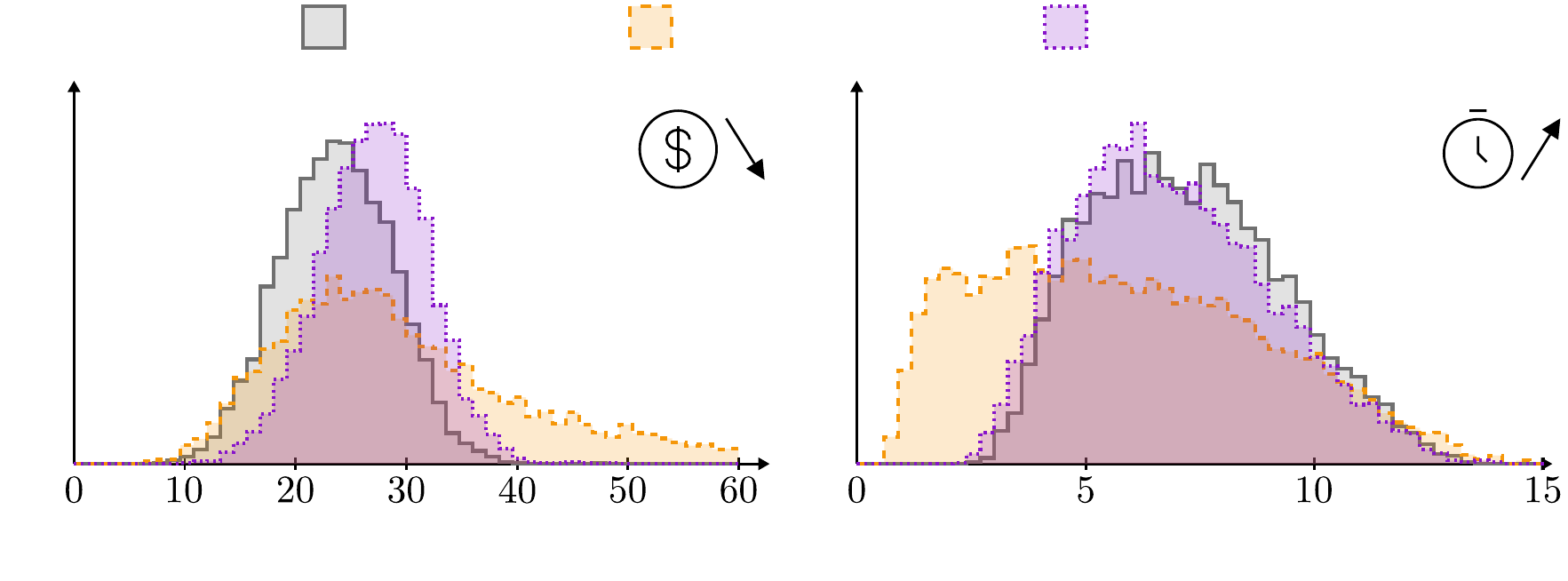}};
    \end{tikzpicture}
    \vspace{-0.8\baselineskip}
    \caption{
        Performance comparison of the three controllers under varying noise distributions $\jointDistributionTrue$, reference distributions $\jointDistributionRef$, and ambiguity radii $\radiusTemplate$.
    }
    \vspace{-.6cm}
    \label{fig:numerical-example}
\end{figure}

%% file: sections/introduction/numerical.tex

\subsection{Why be robust to errors in the noise distributions?}
\label{sec:intro:pendulum}

To motivate the relevance of our \gls*{acr:dr} formulation ahead of its mathematical analysis, we illustrate the benefits of our robust controller using an inverted pendulum with physical parameters taken from \cite{czubarow2025magic}.
This system is commonly used as a canonical benchmark for sensing and control because it is open-loop unstable, lightly damped, and non-minimum phase, which makes it notoriously difficult to control robustly \cite{leong2016understanding}, and challenging even for learning-based controllers \cite{terpin2025usingreinforcementlearningprobe}.

We compare three controllers:
(i) an \textbf{\textcolor{physical}{ideal}} controller with perfect knowledge of the true noise distribution~$\jointDistributionTrue$ affecting the system,
(ii) a \textbf{\textcolor{reference}{nominal}} controller designed for the reference~$\jointDistributionRef$, and
(iii) our \textbf{\textcolor{robust}{robust}} controller, which guards against all distributions in the ambiguity set~$\ambiguitySet$ centered at~$\jointDistributionRef$.
To explore the general trends, we randomly vary $\jointDistributionTrue$, $\jointDistributionRef$, and the radii $\radiusTemplate$, covering
more than 6000
scenarios with low or high noise, accurate or inaccurate references, and conservative or confident settings.
For each configuration\footnote{
    We model $\jointDistributionTrue$ as a zero-mean Gaussian.
    We vary $\jointDistributionTrue$ and $\jointDistributionRef$ by perturbing the entries of the Cholesky factors of their covariances with additive Gaussian noise.
    Radii $\radiusTemplate$ are sampled from a uniform distribution.
    Further details and code are available at
    \url{https://alain.schoebi.ch/dr-lqg}.
},
we then run 500 simulation episodes using the nonlinear model of the pendulum over a 15-second horizon (discretized at 30~Hz).
The histograms in \cref{fig:numerical-example} illustrate how the empirical expected cost (averaged over the balancing period) and the average balancing time (measured until the pendulum falls) vary across the randomized experiments.
While the nominal controller may outperform the robust controller when the reference distribution is, coincidentally, accurate
(\cref{fig:numerical-example}, region~\colorboxSmall{region!10}{\textcolor{region}{R1}}),
the robust controller consistently mitigates \emph{worst-case} scenarios by cutting off the tails observed for the nominal controller (\cref{fig:numerical-example}, region~\colorboxSmall{region!10}{\textcolor{region}{R2}}).
As a result, our robust controller protects against uncertain or misspecified noise~distributions.

%% file: sections/introduction/related.tex

\subsection{Related work}
\label{sec:related}
\begingroup
\newcommand{\nonstationarySymbol}{\dagger}
\newcommand{\ambiguitySetNonstationary}{\ambiguitySet{}{\nonstationarySymbol}}
\newcommand{\policyNonstationary}{\policy^{\nonstationarySymbol}}
\newcommand{\xIndex}{1}
\newcommand{\yIndex}{2}

Our robustness notion is distributional: unlike $\mathcal H_\infty$ or risk-sensitive control, which target energy-bounded disturbances and rare high-cost events, respectively, we hedge against misspecification of the noise distribution via ambiguity sets~\cite{kuhn2025distributionally,tzortzis2016robust,hajar2024distributionally,yan2025distributional,brouillon2025distributionally,li2025distributionally,fochesato2025distributionally,taskesen2024distributionally,taskesen2025optimality}.
In this context, a number of studies~\cite{yan2025distributional,brouillon2025distributionally,li2025distributionally}
restrict to linear or affine policies for tractability without proving their optimality.
Closer to our scope,
\cite{fochesato2025distributionally,taskesen2024distributionally,taskesen2025optimality,cescon2025global}
show that linear policies remain optimal under distributional uncertainty.
While we rely on similar technical tools, these works operate in a more conservative setting where the noise distributions are allowed to vary \emph{over time}, so
that
\mbox{$
    \jointDistribution
    =
    \initialStateDistribution
    \productMeasure
    \processDistribution{0}
    \productMeasure
    \ldots
    \productMeasure
    \processDistribution{\horizon-1}
    \productMeasure
    \;
    \measurementDistribution{0}
    \productMeasure
    \ldots
    \productMeasure
    \measurementDistribution{\horizon-1}
    .
$}
In particular, \cite{taskesen2024distributionally} considers the \gls*{acr:dr} control problem
\begin{equation}
    \label{eq:dr-control-problem:nonstationary}
    \inf{\policy \in \policySet}
    \sup{\jointDistribution \in \ambiguitySetNonstationary}
    \expectedCostPiP
    ,
\end{equation}
where
$
    \ambiguitySetNonstationary
    =
    \ambiguitySetInitialState
    \productMeasure
    \ambiguitySetProcess{0}
    \productMeasure
    \dots
    \productMeasure
    \ambiguitySetProcess{\horizon -1}
    \productMeasure
    \ambiguitySetMeasurement{0}
    \productMeasure
    \dots
    \productMeasure
    \ambiguitySetMeasurement{\horizon -1}
$.

In contrast, our formulation assumes \emph{stationary} disturbances.%
\footnote{
The stationary formulation was first introduced in the preliminary technical notes~\cite{lanzetti2024optimality}, which this letter fully replaces.}
Stationary or slowly time-varying noise models often arise in applications such as chemical processes, robotic systems, and power-electronics systems. In such cases, ignoring stationarity can introduce unnecessary conservatism and degrade performance, justifying a dedicated study of this setting.
The following example compares our approach with~\cite{taskesen2024distributionally}.

\newcommand{\KZeroZero}{\kappa_0}
\newcommand{\KOneZero}{\kappa_1}
\newcommand{\KOneOne}{\kappa_2}
\textbf{Motivating example.}
Consider the linear system
\begin{equation*}
{
\setlength{\arraycolsep}{2pt}
\state{\t+1} =
\begin{bmatrix}
0 & -1 \\
1 & 0
\end{bmatrix}
\state{\t}
+
\begin{bmatrix}
1 \\ 0
\end{bmatrix}
\controlInput{\t} + \processNoise{\t},
\quad
\observation{\t} =
\begin{bmatrix}
0 & 1
\end{bmatrix}
\state{\t} + \measurementNoise{\t}
,
}
\end{equation*}
with horizon $\horizon = 2$, initial state $\state{0} = \initialState = 0$ (taken as known for simplicity of exposition), and cost $10 (\state{2})_{\xIndex}^2 + \controlInput{0}^2 + \controlInput{1}^2$.
Here,~$(\cdot)_{i}$ and $(\cdot)_{ij}$ denote vector and matrix entries.
We restrict to linear policies of the form
$\policy{0}(\observation{0}) = \KZeroZero\observation{0}$
and
$\policy{1}(\observation{0}, \observation{1}) = \KOneOne\observation{1} + \KOneZero\observation{0}$.
The expected cost reduces to
\begin{equation*}
   (
       11 \KOneOne^2
       - 20\KOneOne
       + 10
    ) (\processCovariance{0})_{\yIndex\yIndex}
   +
   11\KOneOne^2 \measurementCovariance{1}
   +
   10(\processCovariance{1})_{\xIndex\xIndex},
\end{equation*}
after noting that $\KZeroZero = \KOneZero = 0$ is always~optimal~(since $\observation{0} = \measurementNoise{0}$ is pure noise).
The matrices $\Sigma_{\bullet}$ denote the covariances of the corresponding noise distributions.
We then construct the ambiguity sets $\ambiguitySet$ and $\ambiguitySetNonstationary$ with standard Gaussians as references and radii $\radiusProcess{\t} = \radiusProcess = 5$, $\radiusMeasurement{\t} = \radiusMeasurement = 1$.

\textit{Stationary setting \eqref{eq:dr-control-problem}.}
Imposing identical noise distributions across time (e.g., $\processCovariance{0} = \processCovariance{1}$) and solving the control problem numerically yields a worst-case expected cost of $372$. This is achieved with $\KOneOne \approx 0.35$, and the worst-case stationary covariances are $\processCovariance \approx \diag{35.6,2.5}$ and $\measurementCovariance = 4$.

\textit{Nonstationary setting \eqref{eq:dr-control-problem:nonstationary}.}
Allowing time-varying noise, we obtain $\KOneOne = \frac{9}{11} \approx 0.82$ with worst-case cost $425$, realized when
$\processCovariance{0} = \diag{1, 6^2}$, $\processCovariance{1} = \diag{6^2, 1}$,
$\measurementCovariance{0} = \measurementCovariance{1} = 4$,
\ie when the disturbances change abruptly across time.

\textit{Discussion.}
As anticipated, the cost in \eqref{eq:dr-control-problem:nonstationary} is higher than in \eqref{eq:dr-control-problem}, since $\ambiguitySet \subseteq \ambiguitySetNonstationary$.
Importantly, the two settings yield different optimal policies, denoted by $\policyStar$ and $\policyNonstationary$;
hence, the stationary formulation is not simply a special case of \cite{taskesen2024distributionally}, but leads to a different optimal controller.
Deploying the nonstationary policy $\policyNonstationary$ under stationary disturbances, or vice versa, results in suboptimal performance.
In particular, using~$\policyNonstationary$ instead of~$\policyStar$ under stationary disturbances strictly increases the worst-case cost from $372$ to $391$.
Thus, ignoring stationarity leads to an overly conservative controller: guarding against implausible time-varying disturbances is unnecessary and degrades performance.
Stationarity is therefore valuable structural information that should be exploited \mbox{when available.}

\endgroup

%% file: sections/theory/import.tex

\section{Theoretical Result}
\label{sec:main-result}

\input{sections/theory/theorem}

\input{sections/theory/discussion}

\input{sections/theory/preliminaries}

\input{sections/theory/proof}

\input{sections/theory/extension/import}

%% file: sections/theory/theorem.tex

\subsection{Optimality of linear policies and Nash equilibrium}
In the theorem below, we show that the \gls*{acr:dr} control problem is solved by
linear policies $\policy$ of the form
$\policy{\t}(\observation{0}, \dots, \observation{\t}) = \sum_{\tTau=0}^\t \controlInputGain{\t}{\tTau} \observation{\tTau}$
and
Gaussian distributions $\jointDistribution$ in the ambiguity set~$\ambiguitySet$.
We denote the corresponding sets by $\policySetLinear$ and $\ambiguitySetGaussian$.

\begin{boxedtheorem}[Optimality of linear policies]
\label{thm:main-result}
When solving the \gls*{acr:dr} control problem~\eqref{eq:dr-control-problem}, it suffices to restrict to the set of linear policies $\policySetLinear$ and Gaussian distributions $\ambiguitySetGaussian$:
\begin{equation*}
    \eqref{eq:dr-control-problem}
    =
    \inf{\policy \in \policySet}
    \sup{\jointDistribution \in \ambiguitySet}
    \expectedCostPiP
    =
    \min{\policy \in \policySetLinear}
    \max{\jointDistribution \in \ambiguitySetGaussian}
    \expectedCostPiP
    .
\end{equation*}
In particular,
the problem admits a Nash equilibrium \mbox{$(\policyStar, \jointDistributionStar)$}
with linear policy $\policyStar \in \policySetLinear$
and Gaussian noise distribution $\jointDistributionStar \in \ambiguitySetGaussian$,
so that for all $\policy \in \policySet$ and $\jointDistribution \in \ambiguitySet$,
\begin{equation*}
    \expectedCost{\policyStar}{\jointDistribution}
    \le
    \expectedCost{\policyStar}{\jointDistributionStar}
    \le
    \expectedCost{\policy}{\jointDistributionStar}
    .
\end{equation*}
Consequently, any Nash equilibrium policy $\policyStar$ solves~\eqref{eq:dr-control-problem}.
\end{boxedtheorem}

%% file: sections/theory/discussion.tex

\textbf{Optimality of linear policies.}
\cref{thm:main-result} states that \emph{linear} policies remain optimal even in this \gls*{acr:dr} setting.
Moreover, the optimal policy $\policyStar$ is a best response to the worst-case distribution $\jointDistributionStar$, which is Gaussian.
Consequently, the optimal controller retains the standard \gls*{acr:lqg} structure, \ie a \gls*{acr:lqr} driven by a \gls*{acr:kf}, both of which can be computed efficiently.

\textbf{Robustness.}
Although both the reference and worst-case distributions are Gaussian, the optimal robust controller remains robust against \emph{all} distributions in the ambiguity set~$\ambiguitySet$, including the non-Gaussian ones (see \cref{fig:ambiguity-set}).

\textbf{Nash vs. Stackelberg games.}
The existence of a Nash equilibrium implies that the zero-sum Stackelberg game is equivalent to a game where the players move simultaneously. In particular, there is no advantage in moving first.
Thus, the Nash equilibrium policy $\policyStar$ is also a Stackelberg equilibrium policy of the sequential game where the engineer moves first, and thereby solves the \gls*{acr:dr} control problem \eqref{eq:dr-control-problem}.

\textbf{Extensions.}
Although we focus on
Wasserstein ambiguity sets, our formulation
extends to the divergence-based ambiguity sets of~\cite{taskesen2025optimality}.
Another extension is to infinite-horizon formulations, where our stationary setting is expected to recover the nonstationary formulation of~\cite[Section 6]{taskesen2025optimality},
because worst-case disturbances happen to be stationary~\cite[Theorem 6.1]{taskesen2025optimality}.

%% file: sections/theory/preliminaries.tex

\subsection{Technical tools}
\label{sec:nash:preliminaries}
Before proving \cref{thm:main-result}, we adopt
the purified-output formulation~\cite{ben2006extending}
to eliminate cyclic dependencies in the inputs.
Using this formulation, we parameterize linear policies \mbox{$\policy \in \policySetLinear$}
via block lower-triangular matrices
\mbox{$\controlInputGainMatrixInR$} as in~\cite{taskesen2024distributionally},
where $\controlInputGainMatrixSet$ denotes the set of all such matrices.
The expected cost~\eqref{eq:the-cost} under a linear policy then simplifies to
$
    \expectedCostPiP
    =
    \expectedCostLinearStationaryUSigmaExpression
    ,
$
where
$
\matBatchCostInitialStateCovariance
,
\matBatchCostProcessCovariance
,
\matBatchCostMeasurementCovariance
\succeq 0
$
depend quadratically on~$\controlInputGainMatrix$ and are defined in
Appendix~\ref{appendix:matrices}.
The matrices $\initialStateCovariance, \processCovariance, \measurementCovariance$ denote the covariances of the marginal distributions of $\jointDistribution$,
which we collect in
\mbox{$\allDisturbancesCovariance = (\initialStateCovariance, \processCovariance, \measurementCovariance)$};
analogously, the reference covariances are collected in
\mbox{$\allDisturbancesCovarianceRef = (\initialStateCovarianceRef, \processCovarianceRef, \measurementCovarianceRef)$}.
Since the cost depends only
on the noise covariances, and not on the full distributions, we write $\expectedCostLinearStationaryUSigma$ to depend directly on
$\allDisturbancesCovariance$ instead of $\jointDistribution$.

Accordingly, we project the ambiguity set~$\ambiguitySet$, defined over~$\jointDistribution$, onto the space of covariance matrices.
This projection yields the covariance ambiguity set
\mbox{$
    \ambiguitySetGelbrichCovariance
    =
    \ambiguitySetGelbrichInitialState
    \times
    \ambiguitySetGelbrichProcess
    \times
    \ambiguitySetGelbrichMeasurement
$},
where
\begin{equation*}
    \ambiguitySetGelbrichTemplate
    =
    \{
        \covarianceTemplate \in \posSemiDef{\textrm{dim}(\templateVariable)}
        :
        \gelbrichDistance{\covarianceTemplate}{\covarianceTemplateRef}
        \le \radiusTemplate
    \},
\end{equation*}
for $\templateVariableInXVW$.
Here $\gelbrichDistance$ denotes the Gelbrich distance between conformable covariance matrices
and is defined by
\begin{equation*}
    \gelbrichDistance{\covarianceMatrixA}{\covarianceMatrixB}
    =
    \sqrt{\trace*{
        \covarianceMatrixA + \covarianceMatrixB -2 \half{\left( \half{\covarianceMatrixB} \covarianceMatrixA \half{\covarianceMatrixB} \right)}
    }}
    .
\end{equation*}
The projection result follows from (i) $\gelbrichDistance$ lower bounding the Wasserstein distance $\wassersteinDistance$ and coinciding with it for Gaussian distributions
\cite[Propositions 2.4 and 3.4]{nguyen2023bridging},
and (ii) the reference distributions being Gaussian.
Moreover,
the projection of the Gaussian ambiguity set $\ambiguitySetGaussian$ also results in $\ambiguitySetGelbrichCovariance$.

%% file: sections/theory/proof.tex

\subsection{Proof of \cref{thm:main-result}}
\label{sec:nash:proof}

\newcommand{\lowerBound}{\diamondsuit}
\newcommand{\upperBound}{\heartsuit}

We borrow the strategy of \cite{taskesen2024distributionally} to prove \cref{thm:main-result}.
Specifically, we bound~\eqref{eq:dr-control-problem} from above and below, and show that the two bounds also coincide in our stationary setting.

\textbf{Upper bound $\upperBound$.} Restricting policies to linear ones gives
\begin{equation*}
\eqref{eq:dr-control-problem}
=
\inf{\policy \in \policySet}
\sup{\jointDistribution \in \ambiguitySet}
\expectedCostPiP
\le
\inf{\policy \in \policySetLinear}
\sup{\jointDistribution \in \ambiguitySet}
\expectedCostPiP
.
\end{equation*}
We use the shorthand $\expectedCostLinearStationary$ and project the ambiguity set~$\ambiguitySet$ onto the covariance space, which yields
\begin{equation*}
\inf{\policy \in \policySetLinear}
\sup{\jointDistribution \in \ambiguitySet}
\expectedCostPiP
=
\inf{\controlInputGainMatrixInU}
\sup{\allDisturbancesCovariance \in \ambiguitySetGelbrichCovariance}
\expectedCostLinearStationaryUSigma
\equalByDefRight
\upperBound
.
\end{equation*}

\textbf{Lower bound $\lowerBound$.} By weak duality, we have
\begin{equation}
\label{eq:security-levels}
\eqref{eq:dr-control-problem}
=
\inf{\policy \in \policySet}
\sup{\jointDistribution \in \ambiguitySet}
\expectedCostPiP
\ge
\sup{\jointDistribution \in \ambiguitySet}
\inf{\policy \in \policySet}
\expectedCostPiP
.
\end{equation}
Restricting the disturbances to Gaussian ones yields
\begin{equation*}
\sup{\jointDistribution \in \ambiguitySet}
\inf{\policy \in \policySet}
\expectedCostPiP
\ge
\sup{\jointDistribution \in \ambiguitySetGaussian}
\inf{\policy \in \policySetLinear}
\expectedCostPiP
,
\end{equation*}
where
we restricted $\policySet$ to $\policySetLinear$ without loss of optimality, since linear policies are optimal under Gaussian disturbances, as per standard \gls*{acr:lqg}.
Similarly, using $\expectedCostLinearStationary$ and projecting the Gaussian ambiguity set $\ambiguitySetGaussian$ onto the covariance space yields
\begin{equation*}
\sup{\jointDistribution \in \ambiguitySetGaussian}
\inf{\policy \in \policySetLinear}
\expectedCostPiP
=
\sup{\allDisturbancesCovariance \in \ambiguitySetGelbrichCovariance}
\inf{\controlInputGainMatrixInU}
\expectedCostLinearStationaryUSigma
\equalByDefRight
\lowerBound
.
\end{equation*}

\textbf{Tightness $\lowerBound=\upperBound$.}
The bounds satisfy $\lowerBound \le \eqref{eq:dr-control-problem} \le \upperBound$.
To prove tightness, we invoke Sion's minimax theorem \cite{sion1958minimax},
which applies since
(i)~$\controlInputGainMatrixSet$ is convex,
(ii)~$\ambiguitySetGelbrichCovariance$ is convex and compact \cite[Lemma A.6]{nguyen2023bridging},
(iii) the cost $\expectedCostLinearStationary$ is jointly continuous, and
(iv)~$\expectedCostLinearStationary$ is linear (hence concave) in~$\allDisturbancesCovariance$ and convex in~$\controlInputGainMatrix$ (as it is nonnegative and quadratic in $\controlInputGainMatrix$).
Therefore, $\lowerBound = \upperBound$.

\textbf{Nash equilibrium.}
Since the bounds coincide, the above chain of inequalities collapses into equalities.
By Berge's maximum theorem and a coercivity argument providing compactness (made explicit in \cref{sec:computational:convergence:decreasing}),
the extrema in~$\lowerBound$ and~$\upperBound$ are attained, and hence so are the outer infimum and supremum in the security levels in~\eqref{eq:security-levels}.
Since the security levels in~\eqref{eq:security-levels} coincide and their outer extrema are attained, a Nash equilibrium exists in $\policySet \times \ambiguitySet$.
Moreover,
by $\lowerBound = \upperBound$,
such an equilibrium is realized in $\policySetLinear \times \ambiguitySetGaussian$.

%% file: sections/theory/extension/import.tex

\subsection{On the zero-mean assumption}

We show that the restriction to zero-mean distributions in
\cref{ass:disturbances} is essential in our setting:
Nash equilibria may fail to exist in the nonzero-mean setting
(or, more generally, when allowing the means of the distributions in the ambiguity set to differ from the reference means).
Since the proof of \cref{thm:main-result} is based on the existence of a Nash equilibrium, the same strategy does not establish optimality
of \emph{affine} policies in the nonzero-mean case.
Importantly, this does not rule out the optimality of affine policies, but shows that establishing such a result would require different arguments.

Accordingly, we drop this restriction and consider the larger
ambiguity set~$\ambiguitySetUnconstrainedMean$, which contains distributions
with nonzero means; we denote by~$\ambiguitySetGaussianUnconstrainedMean$ its Gaussian
counterpart.
In this setting, let $\policySetAffine$ collect affine policies of the form
\mbox{$\policy{\t}(\observation{0}, \dots, \observation{\t}) = \sum_{\tau=0}^\t \controlInputGain{\t}{\tau}\observation{\tau} + \controlInputOffset{\t}$}.
The counterexample below shows that the nonzero-mean analogue
of \cref{thm:main-result} may fail:
a Nash equilibrium in $\policySetAffine \times \ambiguitySetGaussianUnconstrainedMean$ need \emph{not} exist.

\input{sections/theory/extension/counterexample}

\input{sections/theory/extension/discussion}

%% file: sections/theory/extension/counterexample.tex

\textbf{Counterexample.}
Consider the linear system
\begin{equation*}
    \state{\t + 1}
    =
    \state{\t} + \controlInput{\t} + \processNoise{\t}
    ,
    \quad
    \observation{\t}
    =
    \state{\t} + \measurementNoise{\t}
    ,
\end{equation*}
with
horizon $\horizon = 2$,
initial state $\state{0} = \initialState = 0$
(here deterministic for simplicity of exposition),
and cost
$\state{2}^2 + \controlInput{0}^2 + \controlInput{1}^2$.
Let the ambiguity set~$\ambiguitySetGaussianUnconstrainedMean$ be defined by
\mbox{$\processDistributionRef = \normal{0}{1}$},
\mbox{$\radiusProcess = 3$},
and
\mbox{$\measurementDistributionRef = \normal{0}{3^2}$},
\mbox{$\radiusMeasurement = 1$}.
One verifies that the \emph{maximin} problem attains the value
$
    \sup{\jointDistribution \in \ambiguitySetGaussianUnconstrainedMean}
    \inf{\policy \in \policySetAffine}
    \expectedCostPiP
    =
    28
$,
for the Gaussian~$\jointDistributionStar$ with means $\processMean = \measurementMean = 0$ and standard deviations
\mbox{$\sigma_{\processNoise} = \sigma_{\measurementNoise} = 4$}.
The inner minimizer~$\policyStar$ is unique and given by
$\policyStar{0}(\observation{0}) = 0$
and
$\policyStar{1}(\observation{0}, \observation{1}) = -\frac{1}{4}\observation{1}$ (\ie with zero offsets).
If an equilibrium existed, this pair would necessarily constitute one, by uniqueness of~$\policyStar$.
However, $\jointDistributionStar$ is \emph{not} a best response to~$\policyStar$, since the feasible Gaussian with
$\processMean = 3$, $\measurementMean = -1$,
$\sigma_{\processNoise} = 1$, $\sigma_{\measurementNoise} = 3$
yields the strictly larger cost $33.25$.
Hence, no Nash equilibrium exists in $\policySetAffine \times \ambiguitySetGaussianUnconstrainedMean$.

%% file: sections/theory/extension/discussion.tex

\textbf{Discussion.}
This should be contrasted with the related nonstationary Nash equilibrium result of~\cite{taskesen2025optimality}, which allows nonzero-mean noise distributions, but only under the second-moment orthogonality constraint
$\expected{}{z z'^\top} = 0$ for any two distinct time-indexed noise terms.
Under mutual independence, this constraint becomes $\expected{}{z}\expected{}{z'}^\top \!\! = \! 0$, forcing all but possibly one noise marginal to have zero mean.
Thus, in the independent-noise setting considered here, which is standard in control,~\cite{taskesen2025optimality} provides only a limited nonzero-mean extension, and our counterexample shows that, without such a restriction, a Nash equilibrium may fail to exist.
Nevertheless, using~\cite{lanzetti2024variational}, one can show that a best response to affine policies remains Gaussian, suggesting a promising direction for future work.

%% file: sections/computation/import.tex

\section{Computational Method}
\label{sec:computational}

\input{sections/computation/algo}
\input{sections/computation/standalone}
\input{sections/computation/convergence}
\input{sections/computation/equivalence}

%% file: sections/computation/algo.tex

\setlength{\textfloatsep}{10pt}
\algtext*{EndFor}

\begin{algorithm}[t]
\caption{Regularized Iterated Best Response (IBR)}
\label{algo:regularized-ibr}
\begin{algorithmic}
\Require
step sizes~$\stepSizeEngineer{k}, \stepSizeDevil{k} \in \openClosedInterval{0}{1}$, reference covariances $\allDisturbancesCovarianceRef$
\State
Initialize $
    \allDisturbancesCovarianceIterate{0}
    \gets
    \allDisturbancesCovarianceRef
    $,
    $
    \controlInputGainMatrixIterate{0}
    \gets
    0
$

\For{$k = 0, 1, \dots$ \textbf{until convergence}}
    \State \engineerSymbol* \emph{Engineer's turn via LQG}
    \State
    $\displaystyle
        \controlInputGainMatrixBR
        \in
        \colorboxEngineerEquation{\argmin{\controlInputGainMatrixInU}}
        \,
        \expectedCostLinearStationary{\controlInputGainMatrix}{\allDisturbancesCovarianceIterate{k}}
    $
    \State
    $\displaystyle
        \controlInputGainMatrixIterate{k+1}
        \gets
        (1 - \stepSizeEngineer{k}) \controlInputGainMatrixIterate{k}
        +
        \stepSizeEngineer{k}
        \controlInputGainMatrixBR
    $
    \State \devilSymbol* \emph{Devil's turn via bisection}
    \State
    $\displaystyle
        \allDisturbancesCovarianceBR
        \in
        \colorboxDevilEquation{\argmax{\allDisturbancesCovariance \in \ambiguitySetGelbrichCovariance}}
        \,
        \expectedCostLinearStationary{\controlInputGainMatrixIterate{k+\delay}}{\allDisturbancesCovariance}
    $,
    \quad with $\delay \in \{0, 1\}$
    \quad
    \State
    $\displaystyle
        \allDisturbancesCovarianceIterate{k+1}
        \gets
        (1 - \stepSizeDevil{k}) \allDisturbancesCovarianceIterate{k}
        +
        \stepSizeDevil{k} \allDisturbancesCovarianceBR
    $
\EndFor

\noindent\Return final linear policy $\controlInputGainMatrixStar$
\end{algorithmic}
\end{algorithm}

%% file: sections/computation/standalone.tex

We present a computational method to derive an optimal policy $\policyStar$ solving \eqref{eq:dr-control-problem}.
By \cref{thm:main-result}, we can equivalently seek an equilibrium in the game
\eqref{eq:dr-control-problem} with $\policy \in \policySetLinear$ and $\jointDistribution \in \ambiguitySetGaussian$. Under the parametrization of \cref{sec:nash:preliminaries}, this reduces to
\begin{equation}
    \label{eq:expected-cost:with-sets:linear}
    \expectedCostLinearStationaryUSigma
    ,
    \quad
    \text{with}
    \quad
    \controlInputGainMatrix \in \controlInputGainMatrixSet
    ,\quad
    \allDisturbancesCovariance \in \ambiguitySetGelbrichCovariance
    .
\end{equation}
To compute an equilibrium $(\controlInputGainMatrixStar, \allDisturbancesCovarianceStar)$, we use
an \gls*{acr:ibr}-like algorithm in which the engineer and the Devil iteratively update their actions, \colorboxEngineer{\emph{minimizing}} or \colorboxDevil{\emph{maximizing}} the
cost $\expectedCostLinearStationary$.
The players may respond to each other simultaneously or in alternating turns,
which \cref{algo:regularized-ibr} captures via the parameter $\delay \in \{0, 1\}$ determining whether the Devil responds to~$\controlInputGainMatrixIterate{k}$ or~$\controlInputGainMatrixIterate{k+1}$.
If the procedure converges, the iterates yield a Nash equilibrium.
To stabilize the iterations, we use damped updates with step sizes~$\stepSizeEngineer{k}, \stepSizeDevil{k} \in \openClosedInterval{0}{1}$.

\textbf{\engineerSymbol* Engineer's turn.}
Since $\allDisturbancesCovarianceIterate{k}$
characterizes Gaussian noise distributions,
the minimization problem reduces to a standard \gls*{acr:lqg} problem.
Only the \gls*{acr:kf} depends on the noise distribution, so the \gls*{acr:lqr} need not be recomputed.

\textbf{\devilSymbol* Devil's turn.}
Since both $\expectedCostLinearStationary$ and $\ambiguitySetGelbrichCovariance$ are separable,
the maximization decomposes into
$
\max{\covarianceTemplate \in \ambiguitySetGelbrichTemplate} \trace{\matBatchCostTemplateCovariance\covarianceTemplate}
$,
for \mbox{$\templateVariableInXVW$}.
Since the objective is linear in $\covarianceTemplate$ and the set~$\ambiguitySetGelbrichTemplate$ is convex, each subproblem is a convex optimization problem that can be solved efficiently via dualization.
If $\matBatchCostTemplateCovariance \neq 0$,
then by \cite[Proposition A.4(ii)]{nguyen2023bridging}
the optimal covariance~$\covarianceTemplateBR$ is unique and given by
$
\covarianceTemplateBR
=
\dualTemplate^2
(\dualTemplate \identity - \matBatchCostTemplateCovariance)^{-1} \covarianceTemplateRef(\dualTemplate \identity - \matBatchCostTemplateCovariance)^{-1}
$.
Here,
$\dualTemplate > \maxeigenvalue{\matBatchCostTemplateCovariance}$
is found by bisection%
\mbox{\cite[Algorithm 2]{nguyen2023bridging}}
as the unique solution of
$
\trace{
    \covarianceTemplateRef
    \matBatchCostTemplateCovariance^{2}
    (\dualTemplate \identity - \matBatchCostTemplateCovariance)^{-2}
}
=
\radiusTemplate^2
$.
If $\matBatchCostTemplateCovariance = 0$,
the maximizer is
arbitrary; we set it to $\covarianceTemplateBR \gets \covarianceTemplateRef$.
We now establish convergence of \cref{algo:regularized-ibr},
which generally fails under \emph{exact} updates, \ie $\stepSizeEngineer{k} = \stepSizeDevil{k} = 1$, but
can be guaranteed under \emph{damped} best-response dynamics.
\begin{boxedtheorem}[Convergence]
\label{thm:convergence}
Under either step-size and update scheme below,
any accumulation point of the iterates
$(\controlInputGainMatrixIterate{k+\delay},\allDisturbancesCovarianceIterate{k})$
of \cref{algo:regularized-ibr}
is a Nash equilibrium
of the game in~\eqref{eq:expected-cost:with-sets:linear},
hence solving \mbox{the control problem~\eqref{eq:dr-control-problem}.}
\begin{enumerate}[label=(\alph*)]
    \item
    $\stepSizeEngineer{k} = \stepSizeDevil{k}$, $\stepSizeDevil{k} \searrow 0$,
    $\sum_{k=0}^\infty \stepSizeDevil{k} = \infty$,
    with $\delay=0$;\label{step:decreasing}
    \item
    $\stepSizeEngineer{k} = \, 1$,\,\, $\stepSizeDevil{k} \searrow 0$,
    $\sum_{k=0}^\infty \stepSizeDevil{k} = \infty$,
    with $\delay=1$.\label{step:fw}
\end{enumerate}
In particular, the values
$\expectedCostLinearStationary{\controlInputGainMatrixIterate{k+\delta}}{\allDisturbancesCovarianceIterate{k}}$
converge.
\end{boxedtheorem}
\newcommand{\refStepDecreasing}{\textit{\ref{step:decreasing}}\xspace}
\newcommand{\refStepFW}{\textit{\ref{step:fw}}\xspace}

We prove \cref{thm:convergence} for the two regimes in \cref{sec:computational:convergence:decreasing,sec:computational:convergence:fw}.
For \refStepDecreasing, with \emph{simultaneous} updates and decreasing step sizes, convergence follows from an analysis of the \gls*{acr:ibr} dynamics.
For \refStepFW, with \emph{alternating} updates and a \emph{constant} engineer step size equal to $1$ (no damping), convergence follows from an equivalence with a \gls*{acr:fw} algorithm~\cite{frank1956algorithm}.
This connection bridges the gap and unifies
\gls*{acr:ibr}-like
\cite{fochesato2025distributionally}
methods and \gls*{acr:fw} schemes
\cite{taskesen2024distributionally,taskesen2025optimality}
in \gls*{acr:dr}-\gls*{acr:lqg} settings.

%% file: sections/computation/convergence.tex

\subsection{Proof of \cref{thm:convergence}\refStepDecreasing}
\label{sec:computational:convergence:decreasing}

\input{sections/computation/technicalities}

To prove \cref{thm:convergence}\refStepDecreasing, we invoke
\cite[Theorem]{hofbauer2006best}, which establishes convergence for continuous-time best-response dynamics, and, via \cite[Proposition 7]{hofbauer2006best}, for the discrete simultaneous-update case.
These apply since (i) $\expectedCostLinearStationary$ is convex in~$\controlInputGainMatrix$ and concave in~$\allDisturbancesCovariance$, (ii)~$\expectedCostLinearStationary$ is jointly continuous,
and
(iii)~$\controlInputGainMatrixSetCompact$
and~$\ambiguitySetGelbrichCovariancePosDef$ are convex and compact.

%% file: sections/computation/technicalities.tex

The minimax arguments below and the proof of case~\refStepFW require compact action spaces and positive-definite covariance matrices.
In line with this, we introduce the auxiliary action spaces
$\controlInputGainMatrixSetCompact$ and $\ambiguitySetGelbrichCovariancePosDef$,
and show that we can
restrict to those spaces
without loss of generality.
We define
\mbox{$
\ambiguitySetGelbrichTemplatePosDef
=
\{
    \covarianceTemplate
    \in
    \ambiguitySetGelbrichTemplate
    \,|\,
    \covarianceTemplate
    \succeq
    \mineigenvalue{\covarianceTemplateRef}
    \identity
    \succ
    0
\}
$}
and construct
\mbox{
$
\ambiguitySetGelbrichCovariancePosDef =
\ambiguitySetGelbrichInitialStatePosDef
\times
\ambiguitySetGelbrichProcessPosDef
\times
\ambiguitySetGelbrichMeasurementPosDef
$},
which is compact and convex.
We then choose $\controlInputGainMatrixSetCompact$ as the closed convex hull of the set of all best responses $\controlInputGainMatrixBR$
to all \mbox{$\allDisturbancesCovariance \in \ambiguitySetGelbrichCovariancePosDef$}.
The set $\controlInputGainMatrixSetCompact$ is compact since~$\expectedCostLinearStationary$ is uniformly coercive in $\controlInputGainMatrix$ over $\ambiguitySetGelbrichCovariancePosDef$%
\mbox{\cite[Lemma 6]{fochesato2025distributionally}}.%

We now show that any Nash equilibrium
$(\controlInputGainMatrixStar, \allDisturbancesCovarianceStar)$
in
$\controlInputGainMatrixSetCompact \times \ambiguitySetGelbrichCovariancePosDef$
is also an equilibrium of the original game~\eqref{eq:expected-cost:with-sets:linear}.
For this, it suffices to show that the actions $\controlInputGainMatrixStar$ and $\allDisturbancesCovarianceStar$ remain best responses over the original action spaces $\controlInputGainMatrixSet$ and $\ambiguitySetGelbrichCovariance$.
Showing that $\controlInputGainMatrixStar$ is still a best response over $\controlInputGainMatrixSet$ is immediate from the definition of~$\controlInputGainMatrixSetCompact$: $\controlInputGainMatrixSetCompact$ contains all best responses to all $\allDisturbancesCovariance \in \ambiguitySetGelbrichCovariancePosDef$, and in particular those to $\allDisturbancesCovarianceStar$.
Similarly, $\allDisturbancesCovarianceStar$ remains a best response over $\ambiguitySetGelbrichCovariance$: by~\cite[Proposition A.4(iii)]{nguyen2023bridging},
for any fixed~$\controlInputGainMatrixStar$, the maximization over~$\ambiguitySetGelbrichCovariance$ admits an optimizer in~$\ambiguitySetGelbrichCovariancePosDef$.
Hence enlarging $\ambiguitySetGelbrichCovariancePosDef$ to~$\ambiguitySetGelbrichCovariance$ does not improve the maximum value, and $\allDisturbancesCovarianceStar$ remains optimal. %
Finally, the existence of a Nash equilibrium in $\controlInputGainMatrixSetCompact \times \ambiguitySetGelbrichCovariancePosDef$ follows from Sion's minimax theorem, as in \mbox{\cref{sec:nash:proof}}.
We henceforth work with $\controlInputGainMatrixSetCompact$ and~$\ambiguitySetGelbrichCovariancePosDef$.
Importantly, this restriction leaves the engineer's and the Devil's turns unchanged, and therefore does not affect the computations of \cref{algo:regularized-ibr}.

%% file: sections/computation/equivalence.tex

\subsection{Proof of \cref{thm:convergence}\refStepFW}
\label{sec:computational:convergence:fw}

We show that $\stepSizeEngineer{k} = 1$ with alternating updates reduces the algorithm to a \gls*{acr:fw} algorithm applied to the \emph{maximin} problem
\begin{equation*}
    \max{\allDisturbancesCovariance \in \ambiguitySetGelbrichCovariancePosDef}\minCostFunctionSigmas, \quad
    \text{with}\quad
    \minCostFunctionSigmas = \min{\controlInputGainMatrixInUCompact} \expectedCostLinearStationaryUSigma
    .
\end{equation*}
\gls*{acr:fw} thus maximizes $\minCostFunction$ over $\ambiguitySetGelbrichCovariancePosDef$.
To prove that the \gls*{acr:fw} iterates, initialized at
$\allDisturbancesCovarianceIterate{0} \gets \allDisturbancesCovarianceRef$,
coincide with those of \cref{algo:regularized-ibr},
we show that the \gls*{acr:fw} steps match the updates in \cref{algo:regularized-ibr}.

\textbf{Linearization step.}
The \gls*{acr:fw} algorithm first linearizes~$\minCostFunction$ at the current iterate $\allDisturbancesCovarianceIterate{k}$.
Since $\minCostFunction$ is defined as the pointwise minimum of \emph{linear} functions $\expectedCostLinearStationary{\controlInputGainMatrix}{\cdot}$,
the linearization is given by the ``active'' linear function $\expectedCostLinearStationary{\controlInputGainMatrixBR}{\cdot}$,
where
\mbox{%
$
\controlInputGainMatrixBR
\in
\argmin{\controlInputGainMatrixInUCompact}\expectedCostLinearStationary
{\controlInputGainMatrix}
{\allDisturbancesCovarianceIterate{k}}
$
}%
is the corresponding minimizer for~$\allDisturbancesCovarianceIterate{k}$.
This follows from Danskin's theorem~\cite[Proposition A.3.2]{bertsekas2009convex},
using the compactness of~$\controlInputGainMatrixSetCompact$
and uniqueness of~$\controlInputGainMatrixBR$.
The minimizer is unique because
$\expectedCostLinearStationary$
is a positive definite quadratic function of $\controlInputGainMatrix$
under
$\measurementCovariance, \stageInputCostMatrix{\t} \succ 0$.
Accordingly, the \gls*{acr:fw} algorithm computes the best response $\controlInputGainMatrixBR$ to $\allDisturbancesCovarianceIterate{k}$,
corresponding to the engineer's turn in \cref{algo:regularized-ibr},
where the update
with $\stepSizeEngineer{k} = 1$
yields
$\controlInputGainMatrixIterate{k+1} \gets \controlInputGainMatrixBR$.

\textbf{Maximization step.}
The \gls*{acr:fw} algorithm then maximizes the linearized objective over $\ambiguitySetGelbrichCovariancePosDef$, which reduces to
$
    \allDisturbancesCovarianceBR
    \in
    \argmax{
        \allDisturbancesCovariance \in \ambiguitySetGelbrichCovariancePosDef
    }
    \expectedCostLinearStationary{
        \controlInputGainMatrixBR
    }{
        \allDisturbancesCovariance
    }
    $,
and coincides with the Devil's~turn for $\delay = 1$.
The next covariance iterate $\allDisturbancesCovarianceIterate{k+1}$ is then obtained as in \cref{algo:regularized-ibr} via the
step size $\stepSizeDevil{k}$.
Since the two \gls*{acr:fw} steps match the engineer's and the Devil's turns,
the two algorithms produce
identical next iterates
$\allDisturbancesCovarianceIterate{k+1}$ and $\controlInputGainMatrixIterate{k+1}$,
where $\controlInputGainMatrixIterate{k+1}$ in \gls*{acr:fw} denotes the minimizer of the linearization step at~$\allDisturbancesCovarianceIterate{k}$.
Thus, by induction, the two algorithms are equivalent.

\textbf{Convergence.}
Under the step-size conditions stated in~\refStepFW for $\stepSizeDevil{k}$, the \gls*{acr:fw} algorithm converges in value for concave $L$-smooth functions over compact convex sets~\cite[Theorem 3.3]{xu2017convergence}.
Concavity of~$\minCostFunction$ follows from the linearity of~$\expectedCostLinearStationary$ in~$\allDisturbancesCovariance$, and
smoothness can be shown via \cite[Proposition 4.2]{taskesen2024distributionally}.
Hence, any limit point of the iterates is a maximizer~$\allDisturbancesCovarianceStar$ of~$\minCostFunction$.
If~$\controlInputGainMatrixStar$ denotes the corresponding inner minimizer, which is unique, then \mbox{$(\controlInputGainMatrixStar, \allDisturbancesCovarianceStar)$} forms a Nash equilibrium,
\mbox{completing the proof.}

\textbf{Remark.} The suboptimality gap decays at the rate $\mathcal{O}(1/k)$ by \cite[Theorem 4.10 and Corollary 4.8]{xu2017convergence}.
Comparable rates for conditions~\refStepDecreasing in \cref{sec:computational:convergence:decreasing} are known for continuous-time best-response dynamics or in finite-action settings~\cite{hofbauer2006best}.
Although these results do not directly apply, we conjecture that the same rate holds here.
Finally, since the simultaneous and alternating schemes differ only by a one-step delay, we also expect \gls*{acr:ibr} and \gls*{acr:fw} convergence to extend to the complementary schemes via perturbation and stale-linearization arguments, respectively; details are omitted for brevity.

%% file: sections/conclusion.tex

\section{Conclusion}
Linear systems affected by stochastic disturbances are ubiquitous, yet accurately characterizing their noise distributions remains challenging in practice, making distributionally robust control essential.
Exploiting the stationarity of disturbances, we derive a controller that is less conservative than recent approaches.
Beyond its intuitive computational method, our results provide a foundation for several future research directions, including deployment on real-world systems such as~\cite{czubarow2025magic}, investigation of the optimality of affine policies in the nonzero-mean setting, and online adaptation of the ambiguity radii as more data on the noise becomes available.

%% file: sections/appendix/import.tex

\appendix

\input{sections/appendix/matrices}

%% file: sections/appendix/matrices.tex

\section{Convenient matrix definitions}
\label{appendix:matrices}

\newcommand{\Aprod}[2]{\mathcal{A}_{#1}^{#2}}

Batch system matrices \cite{taskesen2024distributionally}:
\vspace{-0.5\baselineskip}
\begin{align*}
    \matBatchC &=
    {\setlength{\arraycolsep}{1pt}
    \begin{bmatrix}
        \Cmat{0} & 0 \\
        & \ddots & \ddots \\
        & & \Cmat{\horizon-1} & 0
    \end{bmatrix}
    }
    ,
    \;\,
    \matBatchG
    =
    {\setlength{\arraycolsep}{2pt}
    \begin{bmatrix}
                         0 &                                    \\
              \Aprod{1}{1} &                                    \\
                    \vdots & \ddots &                           \\
       \Aprod{1}{\horizon} & \cdots & \Aprod{\horizon}{\horizon}
    \end{bmatrix}
    }
    ,
    \;\,
    \matBatchL
    =
    {
    \begin{bmatrix}
               \Aprod{0}{0} \\
                     \vdots \\
        \Aprod{0}{\horizon}
    \end{bmatrix}
    }
    ,
\end{align*}
with $\Aprod{\tTau}{\t} = \Amat{\t-1} \cdots \Amat{\tTau+1}\Amat{\tTau}$ for $\tTau < \t$, and $\Aprod{\tTau}{\t} = \identity$ otherwise,
$
\matBatchQ
=
\diag{\stageStateCostMatrix{0}, \dots, \terminalStateCostMatrix}
$,
$
\matBatchR
=
\diag{\stageInputCostMatrix{0}, \dots, \stageInputCostMatrix{\horizon-1}}
$,
$
\matBatchB
=
\diag{\Bmat{0}, \dots, \Bmat{\horizon-1}}
$,
and
$
\matBatchH
=
\matBatchG \matBatchB
$.

Cost matrices:
\begin{align*}
    \matBatchCostInitialStateCovariance &= \matBatchCostInitialStateCovarianceExpressionCompact
    ,
    \\
    \matBatchCostProcessCovarianceTrajectory &= \matBatchCostProcessCovarianceExpressionCompact
    ,
    \\
    \matBatchCostMeasurementCovarianceTrajectory &= \matBatchCostMeasurementCovarianceExpressionCompact
    ,
    \\
    \matBatchCostProcessCovariance &=
    \textstyle{\sum_{\t = 0}^{\horizon-1}}
    \left(\matBatchCostProcessCovarianceTrajectory\right)_{\t, \t},
    \quad
    \matBatchCostMeasurementCovariance =
    \sum_{\t = 0}^{\horizon-1}
    \left(\matBatchCostMeasurementCovarianceTrajectory\right)_{\t, \t},
\end{align*}
where $(\cdot)_{\t, \t}$ denotes the $t$-th diagonal subblock.